\documentclass[10pt, reqno, a4paper]{article}

\usepackage[utf8]{inputenc}
\usepackage{authblk}
\usepackage{amsfonts, amsthm, amsmath, amssymb}
\usepackage{thmtools}
\usepackage[hypertexnames=false,colorlinks=true,linkcolor=green]{hyperref}
\usepackage{cleveref}
\usepackage{url}

\usepackage{xcolor}
\usepackage[normalem]{ulem}
\usepackage{tikz}
\usetikzlibrary{arrows.meta}
\usepackage{pgfplots}
\pgfplotsset{compat=1.18}
\usepackage{pgfplotstable}
\usepackage{graphicx}
\usepackage[caption=false]{subfig}

\newtheorem{thm}{Theorem}[section]
\newtheorem{cor}[thm]{Corollary}
\newtheorem{lem}[thm]{Lemma}

\newtheorem{defn}[thm]{Definition}
\newtheorem{rem}[thm]{Remark}
\numberwithin{equation}{section}

\newcommand{\eps}{\epsilon}

\newcommand{\oI}{\overline{I}}
\newcommand{\oS}{\overline{S}}
\newcommand{\oT}{\overline{T}}
\newcommand{\oU}{\overline{U}}
\newcommand{\oV}{\overline{V}}
\newcommand{\oX}{\overline{X}}
\newcommand{\oY}{\overline{Y}}
\newcommand{\oa}{\overline{a}}
\newcommand{\ob}{\overline{b}}
\newcommand{\oc}{\overline{c}}
\newcommand{\om}{\overline{m}}
\newcommand{\ogamma}{\overline{\gamma}}

\newcommand{\N}{\ensuremath{\mathbb{N}}}
\newcommand{\R}{\ensuremath{\mathbb{R}}}
\newcommand{\C}{\ensuremath{\mathbb{C}}}

\title{On Milstein-Type Methods for Free Stochastic Differential Equations}

\providecommand{\keywords}[1]
{
	\small	
	\noindent\textbf{\textbf{Keywords }} #1
}

\providecommand{\amscodes}[1]
{
	\small	
	\noindent\textbf{\textbf{AMS Codes }} #1
}
\date{\today}

\author[1, 2]{Georg Schluechtermann}
\author[2]{Michael Wibmer\thanks{Correspondig author email: michael.wibmer@hm.edu}}

\affil[1]{\footnotesize Faculty of Mathematics, Informatics and Statistics, LMU Munich, Germany}
\affil[2]{\footnotesize Faculty of Mechanical, Aeronautical and Automotive Engineering\\ University of Applied Sciences, Munich, Germany}

\begin{document}
	
\maketitle

\begin{abstract}
Previously, the authors derived an analog of the Euler-Maru\-yama method (fEMM) for free stochastic differential equations (fSDEs) and proved strong convergence of order $\gamma=0.5$ in $L_1(\varphi)$-norm under certain assumptions. In this paper, we study the development of numerical methods for fSDEs which show strong convergence of order $\gamma=1$ in $L_\infty(\varphi)$. As a side effect, strong convergence of order $\gamma=0.5$ of fEMM can be extended to $L_p(\varphi)$ for $p\in[1,\infty]$. Utilizing the framework of multiple operator integrals (MOI) we derive a stochastic It\^{o}-Taylor expansion of the solution of the fSDE. It is then possible to identify those free stochastic iterated integrals, which must be discretized in order to obtain strong convergence of order $\gamma=1$. The non-commutativity imposes additional difficulties showing that the iterated free stochastic integrals can be simulated directly only under special situations different from the commutative case. We will show, which diffusion terms lead to a Milstein-type method of order $\gamma=1$. For the cases, where a direct calculation is not possible, we approximate the iterated integrals based on a subdivision of the discretization intervals.
As for fEMM, all proposed methods obey strong convergence of order $\gamma=1$ in $L_p(\varphi),\, 1\leq p\leq \infty$. For all methods developed, we show that the numerical solution is uniformly bounded on finite time intervals.

\end{abstract}


\amscodes{46L53, 46L54, 60H10, 65C30}

\maketitle

\noindent
\section{Introduction}
Free stochastic differential equations (fSDE) emerged up based on the free probability theory, developed by D. Voiculescu in the time period of early 1980s to 1990s. This allowed several researchers to show, that in principle the Doeblin-It\^{o}-calculus and tools from random matrix theory (RMT) can be transferred to these non-commutative setting in an appropriate way, see e.g. \cite{kummererspeicher}, \cite{BianezbMATH01003147}, \cite{Biane1998}, \cite{Biane1998-2}, \cite{BIANESPEICHER2001581},  \cite{anshelevic} and references therein. Nevertheless there are certain differences  which indicate that a merely literal translation of the classical stochastic calculus hits its limits. While the classical stochastic differential equations are driven by at most vector-valued stochastic processes (e.g. Brownian motion or L\'{e}vy processes), here, the state space is an abstract von Neumann algebra with unital, normal, faithful trace. The driving process is a so-called free Brownian motion with values in an abstract finite von Neumann algebra. Thus, we have to encounter the non-commutativity. To get an idea, one should consider the von Neumann algebra of $N\times N$-matrices $M_N(\mathbb C)$ equipped with the normalized trace $\frac{1}{N}\mathbb{E}(\text{tr}(\cdot))$.  
Let $\left(W_1^{(N)},\dots, W_d^{(N)}\right)$ be an $d$-tuple of independent, standard $M_N(\mathbb C)$-valued Brownian motions. Then there is a von Neumann algebra $\mathcal{A}$, a trace $\varphi:\mathcal{A}\rightarrow \C$ and $W_1,\dots,W_d$ freely independent, $\mathcal{A}$-valued, self-adjoint processes on $[0,\infty)$ (called semicircular processes) such that
$$
\varphi\left(P\left(W_{i_1}^{(N)}(t_1),\dots,W_{i_r}^{(N)}(t_r)\right)\right)\rightarrow \varphi\left(P(\left(W_{i_1}(t_1),\dots,W_{i_r}(t_r)\right))\right)
$$
almost sure as $N\rightarrow\infty$ for all $i_1,\dots,i_r\in\{1,\dots,d\}$, $t_1,\dots t_r \geq 0$ and polynomials $P$ in $r$ non-commuting variables. An introduction into free probability and the relation to random matrices can be found e.g. in \cite{voicudykemanica}, \cite{MingoSpeicher2017}, \cite{BianezbMATH01003147}. P.~Biane showed in \cite{BIANESPEICHER2001581}, that a certain scaling of the $N\times N$ Hermitian matrix Brownian motion  converges as $N\rightarrow \infty$ to a so-called free Brownian motion. This motivates the viewpoint, that often formulas of $W_1,\dots,W_d$ in an abstract setting can be studied by considering $W_1^{(N)},\dots, W_n^{(N)}$ and then taking $N\rightarrow \infty$.
It was then P.~Biane and R.~Speicher who developed (\cite{Biane1998-2}, \cite{Biane1998}, \cite{BIANESPEICHER2001581}) the framework of stochastic calculus with respect to free Brownian motion. In their introduction they motivated stochastic free calculus from the matrix level with successively taking limits $N\rightarrow\infty$ (although the paper was developed by a different approach). In \cite[Appendix A]{Nikitop-Ito} the author gives a rigorous matrix-valued stochastic calculus and an It\^{o} formula for $C^2$ scalar functions of Hermitian matrix-valued It\^{o}-processes. \\
To the best knowledge of the authors, free stochastic differential equations appeared first in \cite{kummererspeicher} driven by the idea, to obtain new processes from given ones. P.~Biane and R.~Speicher studied in \cite{BIANESPEICHER2001581} diffusion equations driven by Brownian motion on matrix level and in the limit $N\rightarrow\infty$. The free analog of the central limit theorem (see e.g. \cite{MingoSpeicher2017},  \cite{nourdintaqqu}) then shows, that the free stochastic differential equations are a good approximation and helpful modeling tool for the wide-spread used random matrices. A general study on fSDEs was gained by V. Kargin (\cite{kargin}). He approached fSDEs by the Cauchy transform of the solution, obtaining partial differential equations in terms of the resolvent of the solution of the fSDE. Additionally, a Picard-Lindel\"of-type existence result was given.\\
We consider free stochastic differential equations equipped with a single free Brownian motion $(W_t)_{t\geq 0}$, i.e.
\begin{equation}\label{eq:fSDEinIntro}
	dX_t=a(X_t)dt + \sum_{i=1}^db^i(X_t)dW_tc^i(X_t), 
\end{equation}
where $X_t\in\mathcal{A}$ is a self-adjoint operator and $a,b^i, c^i$ are operator valued, adapted and continuous (in operator norm) functions of $X_t$. Examples are the free analog of the Ornstein-Uhlenbeck process $dX_t=\lambda X_t dt + \sigma dW_t$, the so-called geometric Brownian motion I, $dX_t=\theta X_t dt + \sqrt{X_t}dW_t\sqrt{X_t}$, the Geometric Brownian motion II $dX_t=\theta X_t dt + X_tdW_t + dW_tX_t$ (we refer for details to  \cite{kargin}). Note that formulating the process Geometric Brownian motion II via \eqref{eq:fSDEinIntro} we set $d=2$, $b^1(X_t)=c^2(X_t)=X_t$ and $b^2(X_t)=c^1(X_t)=\text{id}$. For details about the formalism, we refer to \Cref{sec:FreeIto}.\\
The Ornstein-Uhlenbeck free process was also studied in \cite{GAO2006177}.  Free Wishart processes were discussed in \cite{capitaine-donati}. A fSDE appeared in the context of free Jacobi processes in \cite{demni}. In \cite{freeCIR}, the authors stated a free variant of Cox–Ingersoll–Ross (CIR) process. In \cite{dabrowski2016free} fSDEs appeared as a tool to construct free analogs of certain transport maps. Matrix-valued stochastic processes as solutions of a  fSDE was the point of study in \cite{Maecki2019UniversalityCF}. Evolution equations in non-commutative probability can be found in the dissertation \cite{JekelEvolution}. For further works on fSDEs see \cite{capitaine-donati}, \cite{demni},  \cite{Maecki2019UniversalityCF}, and \cite{JekelEvolution}.\\
Similar to the classical case, it is understandable that solutions of fSDEs may not be found explicitly. Hence, numerical methods have come into play to obtain approximation solutions to the underlying fSDE. 
At first, one can think of \eqref{eq:fSDEinIntro} as an equation in an abstract $M_N(\C)$-valued algebra. But by the previously mentioned approximations of the distribution of elements $X\in\mathcal{A}$ by matrix-valued elements $X^{(N)}$, one can consider \eqref{eq:fSDEinIntro} in an abstract von Neumann algebra $\mathcal{A}$. It is then natural to approximate the solutions of free stochastic differential equations by numerical methods on a matrix level. Looking at the classical counterpart, we have the Euler-Maruyama as well as the general Milstein scheme at hand (\cite{kloedenplaten}, \cite{milsteintretyakov}), where the latter acts as prototype for attempting first order methods. One of the major questions concerns the speed of convergence of the numerical iterations and here especially in the strong sense. In \cite{SchlueWib2023} the authors developed a free analog of the Euler-Maruyama scheme (fEMM) converging of order $\gamma=1/2$ in the strong sense in $L_1(\varphi)$ and order one in weak sense. These results were extended to defined stochastic theta methods for fSDEs and studying numerical stability (\cite{NiuWeiYinZhen}).\\
In this paper, we first extend the results on the Euler-Mayurama scheme to general $L_p(\varphi)$-spaces ($1\le p\le \infty$). The reason is, that the free analog of the Burkholder-Gundy inequality is also valid in $L_\infty(\varphi)$ (\cite{Biane1998}). Second, we are seeking for Milstein type methods of order $\gamma=1$ in the strong sense. Similar to the classical case the major ingredient in developing higher order methods is a stochastic  It\^{o}-Taylor 
 expansion of the functions of drift and diffusion terms of the underlying fSDE. Here, the deep result of N.~Azamov, A.~Carey, P.~Dodds and F.~Sukochev (\cite{azamov_carey_dodds_sukochev_2009}) on multiple operator integrals (MOI) hits the scene. For details about MOIs, we also refer to the book \cite{Skripka2019}. The It\^{o}-Taylor approximation and the representation of the derivatives in operator sense enables us to formulate a iterated stochastic It\^{o}-Taylor expansion of the solution of the fSDE in order to identify those terms leading to methods of strong order $\gamma=1$. We show that terms of higher order in the It\^{o}-Taylor expansion are represented using  multiple, iterated operator integrals. As in the commutative case, to gain speed of $\gamma=1$ in the strong sense, such iterated integrals and the MOIs need to discretized. Nevertheless, in the non-commutative setting the situation is therefore more complex than for ordinary stochastic differential equations. \\
 In the commutative case, when the SDE is equipped with a single Brownian motion $B_t\in\R$, the double iterated integrals can be calculated directly as $\int_0^t\int_0^sdB_u=\frac{1}{2}\left(B_t^2 -t\right)$ (\cite{kloedenplaten}). For commutative SDEs with multiple Brownian motion, a direct calculation of the iterated integrals is possible for SDEs with special structure, so-called SDEs with commutative noise \cite{kloedenplaten}. For other types of SDEs, the iterated integrals can only be approximated. We refer to \cite{kloedenplaten}, \cite{HallernRoessler} and the references therein for an overview on approximation methods. In the non-commutative setting, in order to calculate iterated integrals directly, first, the It\^{o} formula in product form (\cite[Theorem 4.1.2]{Biane1998} has to be imposed and second, the MOIs need to be discretized. The It\^{o} formula allows to resolve iterated free stochastic integrals into a product. We will see, that due to non-commutativity, the applicability of the product formula to the terms of the It\^{o}-Taylor expansion  is necessarily limited to the case of a single diffusion term, i.e. $d=1$. But it turns out, that additionally, to be able to apply the free It\^{o} formula, it is required to commute factors in the MOIs. The commuting operation must be handled by certain perturbation formulas (e.g. see \cite{Skripka2019}), but it comes with the costs of the appearance of triple operator integrals representing the second derivative of the diffusion term. These extra terms do not appear in the commutative setting.\\
 Therefore a direct calculation of the iterated integrals is only possible in two cases. First, if the functions of the diffusion terms are affine functions of the unknown $X_t$, i.e. the triple operator integrals vanish. Second, if one restricts to the case of $\mathcal{A}=\mathcal{M}^N(\R)_{sa}$, the triple operator integrals can be calculated directly with the help of the spectral distribution of $X_t$ (resp. of the numerical approximation $\oX_t$.).\\
 A numerical method for the general nonlinear case with multiple diffusion terms ($d>1$), requires the approximation of the free iterated integrals. We do this by resolving the integrals on a finer discretization of step size $\delta t$ on each discretization interval $\Delta t$ of the fSDE, where $\delta t \leq \Delta t^2$. We develop the method (fSM) starting from the iterated It\^{o}-Taylor expansion and discretize the terms necessary to achieve first order convergence. We will show, that the method can be simplified for $d=1$. The approximation of the iterated integrals on the subintervals with length $\delta t$ do not require the discretization of the previously mentioned  triple operator integrals (which in fact, is unknown). Finally we give the proof of strong convergence order of $\gamma=1$. The generality of the method comes with extra computational effort.
 All together, to carry over the Milstein scheme (\cite{milsteintretyakov}, \cite{kloedenplaten}) to the non-commutative case is only possible for $d=1$ and affine diffusion coefficient. In this special case, it is possible to obtain a derivative free, first order method, where the iterated integrals are resolved exactly. For $d=1$ and nonlinear diffusion, the derivatives of operator functions can be calculated from the spectral distribution of the solution in matrix level. This allows for a Milstein type method in the non-commutative setting, which does not appear in the classical commutative case.  Finally we give numerical examples to reproduce the theoretical results numerically on matrix level both for $d=1$ and $d>1$.\\
The paper is organized as follows.  \Cref{sec:prelim-Free-Stoch-Calc} and \Cref{sec:fSDE} contains some preliminaries on free stochastic differential equations. \Cref{sec:freeItoFuncForm} presents the technique on multiple operator integrals, which in fact serves for the Taylor-like expansion. We give some technical lemmas, which are required in the following.  \Cref{sec:fEMMinLp} is mainly a summary of \cite{SchlueWib2023} for the definition of strong convergence and the Euler-Maruyama method. It also states the first step of the It\^{o}-Taylor expansion. \Cref{sec:MilsteinTypeOrder1} then shows the second iteration of the It\^{o}-Taylor expansion and the necessary steps to develop a method of order $\gamma=1$. In \Cref{thm:thm-order1-general} we formulate a condition on the numerical method to be of order $\gamma=1$. \Cref{thm:strong-conv-bound-num-sol} gives a proof to the statement, that on finite time intervals and given any convergent numerical method in the strong sense ($\gamma>0$), the numerical is uniformly bounded on finite time intervals.
\Cref{sec:num-approx-m-terms} deals with the design of methods of order $\gamma=1$. We first turns to those cases, for which the iterated integrals can be resolved directly. It is followed by the construction of methods for the general case $d>1$.
\Cref{sec:numex} is devoted to several examples to show the desired convergence rates numerically.
\section{Preliminaries - Free Stochastic Calculus}
\label{sec:prelim-Free-Stoch-Calc}
Consider a classical probability space $(\Omega, \mathcal{F}, \mu)$ and random variables as measurable functions $X:\Omega \rightarrow \R$. By taking an algebraic viewpoint, one can consider the algebra of random variables and their expectations as a fundamental concept. It allows to generalize the classical, commutative probability to more general settings, e.g. where the random variables are non-commutative. The space $\mathcal{M}^N(\C)=L_\infty\left(\Omega,\mu,\text{M}_N(\C)\right)$ builds up a $*$-algebra with the unit matrix as identity and $\varphi(M)=\frac{1}{N}\mathbb{E}(\text{tr}(M))$ as a trace mapping the identity matrix to $1$. The classical notion of expectation is then replaced by the linear function $\varphi$ on $\mathcal{M}^N(\C)$. \\
Free probability was created by C. Voiculescu in mid 1980's by studying properties of von Neumann algebras. He introduced the notion of freeness, which extends the notion of independent random variables to non-commutative setting. He also discovered, that random matrices satisfy the freeness conditions asymptotically. By the help of non-commutative algebras it is possible to develop non-commutative probability theory. The limits $N\rightarrow\infty$ can be handled properly in algebraic structures and lead to fruitful concepts. It turns out that non-commutative probability theory is realized by using operator algebras such as von Neumann algebras. We refer e.g. to \cite{MingoSpeicher2017},  \cite{anderson_guionnet_zeitouni_2009}, \cite{TaoIntroRMT} and references therein, for setting up non-commutative probability theory and relations to random matrices. To be complete, we give the following general definition (see e.g. \cite{MingoSpeicher2017}). 
\begin{defn}
	A non-commutative probability space is a pair $(\mathcal{A}, \varphi)$, where $\mathcal{A}$ denotes a von Neumann operator algebra and $\varphi:\mathcal{A}\rightarrow \C$ a faithful unital normal trace.
\end{defn}
Since we consider von Neumann algebras with a unital, faithful and normal trace $\varphi:\mathcal{A}\rightarrow\C$, we can introduce for $1\leq p<\infty$ a norm on $\mathcal{A}$ by $\|X\|_p=\varphi(|X|^p)^{\frac{1}{p}}$. The Banach space completion of $\mathcal{A}$ by $\|\cdot\|_p$ is denoted by $L_p(\varphi)$ (see e.g. \cite{PISIER20031459}).
Since the trace is finite we may consider $\mathcal{A}$ as a subset of the predual $L_1(\varphi)$ of the von Neumann algebra $\mathcal{A}=L_\infty(\varphi)$. 
By $\| \cdot \|$ we denote the usual operator norm in $\mathcal{A}$. Let $\mathcal{A}^{sa}=\{a\in\mathcal{A}, a^*=a\}$. For a non-commutative random variable $X\in\mathcal{A}^{sa}$,	there is a unique probability measure on $\R$ with compact support having the same moments as $X$ (e.g. \cite{MingoSpeicher2017}, \cite{TaoIntroRMT}). This probability measure is the distribution of the non-commutative random variable X.\\
The non-commutative analog of independence of classical random variables is the concept of free independence, or shortly freeness, of subalgebras of $\mathcal{A}$ (\cite{MingoSpeicher2017}). Let $\mathcal{A}_1,\dots \mathcal{A}_n$ be a family of $n\in\N$ subalgebras of $\mathcal{A}$. They are called freely independent (or simply free) in the sense of Voiculescu, if $\varphi\left(X_1X_2\dots X_m\right)=0$ 
whenever the following conditions
\begin{enumerate}
	\item $X_j\in\mathcal{A}_{i(j)} $, where $i(1)\neq i(2), i(2)\neq i(3), \dots , i(n-1)\neq i(n)$, $j=1,\dots,m$
	\item $\varphi(X_i)=0$ for all $i=1,\dots,n$
\end{enumerate} 
hold (\cite[Definition 11]{MingoSpeicher2017}). If $X\in\mathcal{A}$ is a self-adjoint element, then there is a unique spectral measure $\mu$ on $\R$ so that the moments of $X$ are the same as the moments of the probability measure $\mu$ defined by
$
\varphi(X^k)=\int_\R x^k d\mu(x), 
$
see \cite[pp. 51]{MingoSpeicher2017}. An important role plays the Cauchy transform $G_X$ of $\mu$ defined by
$
G_X(z)=\int_{\R} \frac{d\mu(x)}{x-z},
$
which is an analytic function defined on $\C^+$ with values in $\C^+$. The Cauchy transform $G_X$ is the expectation of the resolvent of $X$, i.e.
$
G_X(z)=\varphi\left(\left(X-z\right)^{-1}\right)
$.\\
The Cauchy transform carries all the properties of the spectral probability distribution of the self-adjoint operator $X$. 
In \cite{BIANESPEICHER2001581} the authors used the Hilbert transform of the solution of the underlying fSDE to obtain detailed information about the distribution of the solution. In \cite{kargin} these results were extended to general fSDEs. They can be handled by a corresponding deterministic partial differential equations of the Cauchy transform $G_X$. We will strongly depend on these results since it allows us to check the numerical results.
\subsection{Free Brownian Motion}
\label{sec:FreeBrownMotion}
Motivated by the concept of classical Brownian motion the definition within non-commutative probability is as follows. Consider a von Neumann algebra $\mathcal{A}$ with a faithful normal trace  $\varphi:\mathcal{A}\rightarrow \C$.
A filtration $\mathbb{F}=(\mathcal{A}_t)_{t \geq 0}$ is a family of von Neumann subalgebras $\mathcal{A}_t$ of $\mathcal{A}$ with $\mathcal{A}_s \subset \mathcal{A}_t$ for $s\leq t$.
A family of elements $(X_t)_{t\geq 0}\subseteq \mathcal{A}$ is called adapted to the filtration $\mathbb{F}$, if $X_t\in\mathcal{A}_t$ for all $t\geq 0$. 
\begin{defn}
	A free Brownian motion $(W_t)_{t\geq 0}$ is a family elements of $\mathcal{A}$ adapted to the filtration $\mathbb{F}=(\mathcal{A}_t)_{t\geq 0}$, which admits the following properties:
	\begin{enumerate}
		\item $W_t$ is a self-adjoint element of $\mathcal{A}$ with semi-circular distribution of mean zero and variance $t$, 		
		\item for all $s,t$ with $s\leq t$, the element $W_t-W_s$ is free of $\mathcal{A}_s$ and has a semi-circular distribution with mean $0$ and variance $t-s$.
	\end{enumerate}
\end{defn}
On finite time intervals $J\subset \R$ the elements $W_t$ are uniformly bounded, i.e. $\sup_{t\in J}\|W_t\|<\infty$.
\begin{rem}
	Above definition is taken from \cite{Biane1998}, assuming a filtered probability space.  The definitions in \cite{nourdin} and \cite{kargin} are synonymous.
\end{rem}
\subsection{Stochastic Integration with Respect to Free Brownian Motion, Free It\^{o} Formula}
\label{sec:FreeStochInt}
So called "free stochastic integrals", i.e. stochastic integrals with respect to free Brownian motion were introduced in \cite{kummererspeicher}, \cite{Biane1998}, \cite{anshelevic}. Due to non-commutativity the stochastic integral is build up by integrands, where operators are multiplied both on the left and right of the integration variable. Such integrals are constructed by first introducing piecewise constant processes as integrands, so called simple biprocesses. By defining an appropriate norm, the vector space of simple biprocesses can be completed to the general space of biprocesses. We briefly summarize the construction of free integrals, for details we refer to \cite[Chapters 2, 3]{Biane1998}. \\
Consider the opposite algebra $\mathcal{A}^{op}$ to $\mathcal{A}$ and a decomposition $0=t_0<t_1<\dots t_m<\infty$. A simple adapted biprocess $U_t$ is a piecewise constant map
$$
U_t=
\begin{cases}
	A_k\otimes B_k,&t_k\leq t < t_{k+1}\\
	0,&t_n\leq t
\end{cases},
$$
where $A_k,B_k\in\mathcal{A}_{t_k}$ according to the filtration $\mathbb{F}$. Following  
\cite[Definition 2.2.1]{Biane1998},  the stochastic integral of $U$ with respect to the Brownian motion $(W_t)_{t\geq 0}$ is the integral
($\Delta W_k=W_{t_{k+1}}-W_{t_k}$)
$$
\int_0^{\infty} U_s \sharp d W_s =\int_0^\infty A_sdW_sB_s := \sum_{k=0}^{m-1} U_{t_k} \sharp\left(\Delta W_k\right)=\sum_{k=0}^{m-1} A_{t_k}^j\left(\Delta W_k\right) B_{t_k}^j.
$$
By bilinearity, one can expand this definition to finite sums $U_t=\sum_{i=1}^n A^i_t\otimes B^i_t$.
An It\^{o}-isometry can be obtained for all adapted simple biprocesses. E.g. for $U=A\otimes B$ and $V=C\otimes D$ the isometry reads
$$
\varphi\left[\int U_t \sharp d S_t \cdot\left(\int V_t \sharp d S_t\right)^*\right]=\int\left\langle U_t, V_t\right\rangle_{L_2(\varphi) \otimes L_2(\varphi)} dt = \int \varphi(A_tC_t^*)\varphi(B_tD_t^*)dt .
$$
The vector space of simple biprocesses can be equipped with the norms
$$
 \|U\|_{\mathcal{B}_p}:=\left(\int\left\|U_t\right\|_{L^p\left(\varphi \otimes \varphi^{o p}\right)}^2 dt\right)^{1 / 2}.
$$
In the case of a selfadjoint biprocess $U=A\otimes B$ the norm is $$
\|A\otimes B\|_{\mathcal{B}_p}= \left(\int\varphi(A_t^2)\varphi(B_t^2)dt\right)^{1/2}.
$$
This vector space can be completed to the space $\mathcal{B}_p$ (resp. for $\mathcal{B}_p^a$ for the closed subspace of adapted processes). This allows to extend the mapping $U\mapsto \int U_t\#dW_s$ isometrically to $\mathcal{B}_2^a\rightarrow L_2( \varphi)$.
Free calculus has the essential property, that free stochastic integrals are bounded operators even for $p=\infty$. In the following, we proceed adapted to our needs. Let $A_t, B_t\in\mathcal{A}$ and $(A_t)_{t\geq 0}, (B_t)_{t\geq 0}$ be continuous  processes adapted to the filtration $\mathbb{F}$. We use the simple notation of the free integrals as 
$
  	\int_{o}^{t} A_s dW_s B_s.
$ 
As stated above, in non-commutative setting it is possible to set up a Burkholder-Gundy martingale inequality even for $p=\infty$ (\cite[Section 3.2.]{Biane1998}):
\begin{equation}\label{ineq:BurkholdGundy}
	\left\| \int_0^t U_s\#dW_s\right\| \leq 2\sqrt{2}\left( \int_0^t \|A_s\|^2 \|B_s\|^2ds \right)^{\frac{1}{2}}.
\end{equation}
Hence, the inequality implies
$\left\| \int_t^{t+\Delta t} A_sdW_sB_s \right\|=O(\sqrt{\Delta t})$.

\subsection{Free It\^{o}-Process}

\label{sec:FreeIto}
\begin{defn}\label{def:FreeItoProcess}
	Let $(W_t)_{t\geq0}$ be a free Brownian motion and $X_0\in \mathcal{A}^{sa}_0$.
	Let $a,b^i,c^i:[0,T]\rightarrow\mathcal{A}^{sa}$ 	
	continuous functions in the operator norm and adapted. 
	A free It\^{o}-process is a self-adjoint, adapted process $(X_t)_{t\geq 0}$ of the form  
	\begin{equation}\label{intro-def-freeIto}
		X_t=X_0 + \int_0^t a(s)ds + \sum\limits_{i=1}^d\int_0^tb^i(s)dW_sc^i(s).
	\end{equation}
\end{defn}
Note, that the product $b^i(s)dW_sc^i(s)$ is not necessarily self adjoint. We require $X_t\in\mathcal{A}^{sa}$, so $b^i, c^i$ cannot be chosen arbitrarily.
It would be natural to define the diffusion terms symmetric, i.e. $\sum_{i=1}^d b_t^idW_tb_t^i$. Considering the diffusion $\int_0^tX_s dW_s + dW_sX_s$, it would be possible to write $X_s dW_s + dW_sX_s = (X_s+\text{id})dW_s(X_s+\text{id}) + (X_s)(-dW_s)(X_s) - dW_s$. This requires mixed positive and negative signs before the free Brownian motion at different summands and would make the formalism more complicated. If one summand shows $b^i\neq c^i$, there must be another summand $c^idW_s b^i$ in order to meet the self adjoint requirement. \\
To keep it simple, we continue with the formalism \eqref{intro-def-freeIto}. For example, to express $\int_0^tX_s dW_s + dW_sX_s$, we set  $d=2$ and $b^1(X_t)=c^2(X_t)=X_t$ and $b^2(X_t)=c^1(X_t)=\text{id}$.

\section{Free Stochastic Differential Equations (fSDEs)}
\label{sec:fSDE}
\begin{defn} \label{def-freeItoProcess}
	Let $0<T\leq\infty$ and $I=[0,T[$.
	$X_0$ denotes a self-adjoint element  in $\mathcal{A}^{sa}$ and $a,b^i,c^i:\mathcal{A}\rightarrow \mathcal{A}$ continuous functions in the operator norm. We call 
	\begin{equation}\label{intro-freeSDE-diffform}
		dX_t=a(X_t)dt+ \sum\limits_{i=1}^d b^i(X_t)dW_tc^i(X_t)
	\end{equation}
	a (formal) free stochastic differential equation (fSDE).	
	A solution of \eqref{intro-freeSDE-diffform} with initial condition $X(0)=X_0$ is a self-adjoint continuous adapted processes $(X_t)_{t\in I}$ with the following properties:
	\begin{enumerate}
		\item $X(0)=X_0$ is a self-adjoint element in $\mathcal{A}_0^{sa}$
		\item $X_t\in\mathcal{A}_t^{sa}$ for all $t\in I$
		\item The equation
		\begin{equation}\label{intro-freeSDE-inform}
			X_t=X_0 + \int_0^t a(X_s)ds + \sum\limits_{i=1}^d\int_0^tb^i(X_s)dW_sc^i(X_s)
		\end{equation} 
		is fulfilled for all $t\in I$.
	\end{enumerate}
\end{defn}

\begin{defn}
	We call a function $f:\R\rightarrow\R$ locally operator Lipschitz, if it is a locally bounded, measurable function such that for all $A>0$, there is a constant $L_f(A)>0$ such that
	\begin{equation}\label{eq:L2Lipschitz-a-Lemma}
		\left\| f(X)-f(Y)\right\| \leq L_f(A)\|X-Y\|,
	\end{equation} for elements $X,Y\in\mathcal{A}^{sa}$ and $\|X\|,\|Y\|<A$. If the constant $L_f$ in  \eqref{eq:L2Lipschitz-a-Lemma} does not depend on $A$, we call $f$ (globally) operator Lipschitz or short operator Lipschitz.
\end{defn}
From above definition it follows immediately, that locally operator Lipschitz functions are continuous.
\begin{rem} Since $(X_t)_{t\in I}$ as adapted to $(\mathcal{A}_t)_{t\in I}$, so is the image $f(X_t)$ for continuous $f:\mathcal{A}\rightarrow\mathcal{A}$ (in operator norm).
\end{rem}
We take the following local existence and uniqueness results from \cite{kargin}.
\begin{thm}
Suppose that $a_i, b_i$, and $c_i$ are locally operator Lipschitz functions and $\bar{X}\in\mathcal{A}^{sa}$. Then, there exist $0<T<\infty$ and a family of operators $(X_t)_{t\in [0,T[}$ uniformly bounded in operator norm, such that $X_0=\bar{X}$, and $(X_t)_{t\in[0,T[}$ is a unique solution of \eqref{intro-freeSDE-diffform}.
\end{thm}
\begin{rem}\label{rem:fSDE-sol-theory}
	The existence proof in \cite{kargin}, originally formulated in operator norm, can easily be formulated in $L_p(\varphi)$,
	due to the validity of the free Burkholder-Gundy inequality for $L_p(\varphi), \, 1\leq p \leq \infty$. The solution $(X_t)_{t\in I}$ is therefore uniformly bounded in  $L_p(\varphi), \, 1\leq p \leq \infty$, which is a significant difference to commutative SDEs, for which boundedness in operator norm is not necessarily given. The boundedness property of $(X_t)_{t\in I}$ on finite time intervals will play a major role in the proofs of strong convergence properties in the following.
\end{rem}
As an initial example consider the free analog of the Ornstein-Uhlenbeck process (see \cite{kargin}) defined by the fSDE
\begin{equation}\label{ornsteinuhlenbeck}
	dX_t=\theta X_tdt + \sigma dW_t, \, t\geq 0, \,\,\theta, \sigma\in\R.
\end{equation}
Spectral information about the solution can be obtained by taking the Cauchy transform $G$ of the self-adjoint element $X_t$.
$G$ fulfills a deterministic partial differential equation (\cite[Proposition 3.7]{kargin}). 
Applying the Stieltjes inversion formula (see \cite{kargin}) to its solution, it is possible to recover the spectral distribution of $X_t$. In the case $\theta<0$ it turns out that the density of $X_t$ is a semicircle distribution with radius 
$$R(t)=\sqrt{\frac{2\sigma^2}{|\theta|}(1-e^{-2|\theta| t})}.$$
For $t\rightarrow \infty$ the probability distribution function (PDF) converges to a semicircle with radius $\sigma \sqrt{\frac{2}{|\theta|}}$. The case $\theta\geq0$ is treated in the same way.

For more examples and discussion about the properties of the spectral distribution of the solution we refer to \cite{kargin}. Important fSDEs with a single diffusion term are
\begin{itemize}
	\item Geometric Brownian Motion: 
	$$
		dX_t=\theta X_tdt + \sqrt{X_t}dW_t\sqrt{X_t}
	$$
	By applying $\varphi$ it follows $\varphi(X_t)=e^{\theta X_t}$. 
	\item As shown in \cite{kargin}, the solution of the following equation explodes in finite time
	$$
		dX_t=kX_tdW_tX_t.
	$$
	If $X_0=I$, then the spectral distribution of $X_t$ is defined for all $t\leq 1/k^2$.
\end{itemize}
Equations with a single diffusion term
\begin{itemize}
	\item Geometric Brownian Motion 2:
	$$
	dX_t=\theta X_tdt + X_tdW_t + dW_tX_t
	$$
	This is an example of as fSDE with two diffusion terms. The Brownian motion is the same in both summands.
	\item Free CIR Equation from \cite{freeCIR},
	$$
	dX_t=(a-bX_t)dt + \sigma/2(\sqrt{X_t}dW_t + dW_t\sqrt{X_t}).
	$$
\end{itemize}
~\\~
\section{Operator Integrals, Taylor Series, It\^{o}-Formula}
In this section we introduce some well-known but necessary ingredients. We start with
\subsection{Operator Integrals}
An important ingredient in the following is the notion of a multiple operator integral, originating in the study of perturbations of operator functions. For functions with special properties, it is possible to express the difference $f(X)-f(Y)$, where $X,Y\in \mathcal{A}$, as a so-called double operator integral. This idea can be extended, such that under certain conditions, the Fr\'{e}chet-derivatives of operator functions can be formulated as operator integrals. We refer to \cite{Skripka2019} and  \cite{azamov_carey_dodds_sukochev_2009}. \\
 The definition of a double operator integrals $T_{f^{[1]}}$, resp. a triple operator integral $T_{f^{[2]}}$ is given in
\cite[Definition 4.1]{azamov_carey_dodds_sukochev_2009} and \cite[Lemma 4.5]{azamov_carey_dodds_sukochev_2009}. Recalling: For $X\in\mathcal{A}^{sa},\,Y\in\mathcal{A}$,
\begin{equation}\label{def:optintgeneral}
	\begin{aligned}
		T_{f^{[1]}}^{X,X}(Y) &= \int_{\Pi^{(2)}}e^{i(s_0-s_1)X} Y e^{is_1X}d\nu_f^{(2)}(s_0,s_1), \\
		T_{f^{[2]}}^{X,X,X}(Y,Y) &=
		\int_{\Pi^{(3)}} e^{i(s_0-s_1)X} Y e^{i(s_1-s_2)X} Y e^{is_2X} d\nu_f^{(3)}(s_0,s_1,s_2).
	\end{aligned}
\end{equation}
The definition of the set $\Pi^{(n)}$ and the measure $\nu_f^{(n)}$ can be found in \cite[Lemma 2.1]{azamov_carey_dodds_sukochev_2009}.\\
In our case, we consider functions, taken from the so-called Wiener space $W_n(\R)=\{f\in C^n(\R):~f^{(k)}, \mathcal{F}f^{(k)}\in L_1(\R), \, k=0,\dots,n\}$, $n\in\N$. If $f\in W_1(\R)$ the double operator integral is well defined (resp. $f\in W_2(\R)$ for the triple operator integrals.) As stated in \cite{Biane1998-2}, $C^2$ functions are locally operator Lipschitz. Since $W_n(\R)\subset C^2(\R)$ for $n\geq 3$, functions $f\in W_n(\R), n\geq 3$ are locally operator Lipschitz in operator norm.
P.~Biane and R.~Speicher (\cite{Biane1998}) developed free calculus with respect to free Brownian motion. \\
For functions $f\in W_n(\R)$, it is possible to give a Taylor approximation with appropriate remainder term \cite[Chapter 5.4]{Skripka2019}, \cite[Corollary 5.8]{azamov_carey_dodds_sukochev_2009}.
For example, let $f \in W_3(\mathbb{R})$. Applying the Taylor series expansion \cite[Corollary 5.8]{azamov_carey_dodds_sukochev_2009}, then the estimation of the remainder \cite[Theorem 5.4.4]{Skripka2019} yields
\begin{equation}\label{eq:uuu1}
	f(X_{i+1})-f(X_i) = T_{f^{[1]}}^{X_i,X_i}(\Delta X)   + T_{f^{[2]}}^{X_i,X_i,X_i}(\Delta X,\Delta X) + O(\|\Delta X\|^3) 
\end{equation}
where $\Delta X=X_{i+1}-X_i$ and $X_{i+1},X_i\in\mathcal{A}$.\\
The following  lemmas are needed in the following. Although they are well known in the literature, we state these lemmas for better readability. \Cref{lem:referee} states that e.g. if a left factor $X\in\mathcal{A}$ is pushed into the operator integral, then due to non-commutativity, a triple operator integral appears. \Cref{lem:discrete-opt-int} allows to replace a double operator integral by a difference of the underlying function, together with an estimation of the reminder. It allows to discretize the double operator integral to be implemented in a numerical method.  \Cref{lem:estim-optint-diff} and \Cref{lem:discstochint_lin_growth} are needed e.g. in the proof of \Cref{thm:strong-conv-bound-num-sol}.
\begin{lem}\label[lemma]{lem:referee} Let $A,X,Y \in \mathcal{A}^{sa}$. If $f\in W_3(\R)$ then 
	\begin{eqnarray*}
		XT_{f^{[1]}}^{A,A}(Y)&=&T_{f^{[1]}}^{A,A}(XY)+T_{f^{[2]}}^{A,A,A}(XA-AX,Y)\\
		T_{f^{[1]}}^{A,A}(Y)Z&=&T_{f^{[1]}}^{A,A}(YZ)+T_{f^{[2]}}^{A,A,A}(Y,AZ-ZA).
	\end{eqnarray*}
\end{lem}
\begin{proof} See \cite[Lemma 4.4.4 and 4.4.5]{Skripka2019} and proof of \cite[Theorem 4.4.8]{Skripka2019}, e.g. \cite[(4.4.13)]{Skripka2019}.
\end{proof}

The following lemma is a simplified version of \cite[Theorem 5.4.4, Theorem 5.4.5]{Skripka2019} resp. \cite[Corollary 5.8]{azamov_carey_dodds_sukochev_2009}. For our purposes, it is sufficient to consider functions  $f\in W_n(\R)$. Since $W_n(\R)\subset B_{\infty 1}^n$ the assumptions of \cite[Theorem 5.4.4]{Skripka2019}  are fulfilled (also for the case $p=\infty$). For the definition of the Besov space $B_{\infty 1}^n$ we refer to \cite[Chapter 2.1]{Skripka2019}.
\begin{lem}\label{lem:discrete-opt-int}
	Let $A, V\in\mathcal{A}^{sa}$,  $f\in W_3(\R)$. Then 
	\begin{equation}\label{eq:discrete-opt-int}
		T_{f^{[1]}}^{A,A}(V) = f(A+V) - f(A) - R_{2,f,A}(V)
	\end{equation}
	where $\|R_{2,f,A}(V)\| = \mathcal{O}(\|V\|^2)$. $R_{2,f,A}(V)$ is the Taylor remainder as defined in \cite[(5.4.7)]{Skripka2019}. 
\end{lem}
Note that we rely mainly on the norm estimation of the remainder $R_{2,f,A}(V)$ but not on the specific expression.

\begin{lem}\label{lem:estim-optint-diff}
	Let $A, B, V\in \mathcal{A}^{sa}$ and $f\in W_3(\R)$. Then there is a constant $K_f>0$, such that
	\begin{equation}\label{lem:difference-operator-int}
		\|T_{f^{[1]}}^{A,A}(V)-T_{f^{[1]}}^{B,B}(V)\| \leq K_f \|V\|\|A-B\|.
	\end{equation}
\end{lem}
\begin{proof}
	The perturbation formula \cite[Theorem 4.3.14]{Skripka2019} gives the relation
	$$
	T_{f^{[1]}}^{A,A}(V)-T_{f^{[1]}}^{B,B}(V) = 
	T_{f^{[2]}}^{A,B,A}(A-B,V)-T_{f^{[2]}}^{B,A,B}(V,A-B).
	$$
	The estimation of the triple operator integrals via \cite[Theorem 4.3.8]{Skripka2019} yields the assertion.
\end{proof}

\begin{lem}\label{lem:estim-optint-diff-triple}
	Let $A, B, S, U,\oU,\oS \in \mathcal{A}^{sa}$ and $f\in W_4(\R)$. Let $S,U,\oU,\oS,A,B$ be uniformly bounded in operator norm. Then there is a constant $K>0$, such that
	\begin{equation*}\label{lem:difference-operator-int-2nd}
		\|T_{f^{[2]}}^{A,A,A}(S,U)-T_{f^{[2]}}^{B,B,B}(\oS,\oU)\| \leq 
		K_1 \|S-\oS\|+ K_2\|A-B\| + K_2\|U-\oU\|
	\end{equation*}
\end{lem}
\begin{proof}
	\begin{align*}
		&T_{f^{[2]}}^{A,A,A}(S,U)-T_{f^{[2]}}^{B,B,B}(\oS,\oU) =\\
		&=T_{f^{[2]}}^{A,A,A}(S,U)-T_{f^{[2]}}^{A,A,A}(\oS,U) + T_{f^{[2]}}^{A,A,A}(\oS,U) -T_{f^{[2]}}^{B,B,B}(\oS,\oU) \\
		&=T_{f^{[2]}}^{A,A,A}(S-\oS,U)+ T_{f^{[2]}}^{A,A,A}(\oS,U) -T_{f^{[2]}}^{B,B,B}(\oS,\oU)  \\
		&=T_{f^{[2]}}^{A,A,A}(S-\oS,U)+T_{f^{[3]}}^{A,B,A,A}(A-B,\oS,U) +\\
		&+T_{f^{[3]}}^{B,A,B,A}(\oS,A-B,U)+T_{f^{[3]}}^{B,B,A,B}(\oS,U,A-B)+T_{f^{[2]}}^{B,B,B}(\oS,U-\oU)
	\end{align*}
	Since operator integrals are bounded by their arguments the statement follows. The estimation of the triple operator integrals via \cite[Theorem 4.3.8]{Skripka2019} yields the assertion.
\end{proof}

\begin{lem}\label{lem:discstochint_lin_growth} Let $b,c : \mathcal{A}^{sa}\rightarrow\mathcal{A}^{sa}$ be two locally operator Lipschitz functions . If $X\in\mathcal{A}^{sa}$ is free from the increment of a Brownian motion $\Delta W=W_{t+\Delta t}-W_t$, then there is a constant $L_{bc}>0$, such that the estimation 
	$$\|b(X)\Delta W c(X)\|^2\leq 8L_{bc}(1+\|X\|^2)^2\Delta t$$
	holds.
\end{lem}
\begin{proof}
	We shorten $b=b(X)$, analog for $c$. Considering that the product $b \Delta W c$ is self-adjoint, further applying \cite[Lemma 3.3]{kargin} and using free Burkholder-Gundy
	\begin{equation}\label{eq:proof-bound-num-gronwall-1}
		\|b\Delta W_t c\|^2
		= \left\|\int_{\Delta t} b dW_s c\right\|^2 \le 8 \int_{\Delta t} \|b\|^2\|c\|^2ds \leq 
		8L_{b}L_c\left(1+\|X\|^2\right)^2\Delta t.
	\end{equation}
\end{proof}

\subsection{Free It\^{o} Formula }\label{sec:freeItoFuncForm}
An important ingredient in the development of numerical methods for fSDEs and their convergence properties is a free analog of the It\^{o} formula. Influenced by the formalism from the commutative Wagner-Platen extension for solutions of SDEs (\cite{kloedenplaten}), we  adapt the notation of the It\^{o}-formula introduced in \cite[Proposition 4.3.4]{Biane1998}.
\begin{thm}[Free It\^{o} Formula in Integral Form]\label{freeItoTheorem} Suppose $a, b^i, c^i: \mathcal{A}\rightarrow\mathcal{A}$ are continuous functions  in operator norm for $i=1,\dots,d$. Furthermore $b^i,c^i$ are so that the sum $\sum_{i=1}^db^i(X_t)dW_tc^i(X_t)$ is self-adjoint (see also \Cref{sec:FreeIto}). Let $(X_t)_{t\geq 0 }$ be an adapted, continuous, self-adjoint and free It\^{o}-process and $X_0 \in \mathcal{A}^{sa}_0$. Then for functions $f\in W_3(\R)$, it follows that
	\begin{equation}\label{eq:freeItoFormula}
		f(X_t)=f(X_0)+\int_0^t L^0\left[f(X_s)\right]ds + L^1\left[f(X_s)\right]_0^t
	\end{equation}
	where the operators $L^0,L^1:\mathcal{A}^{sa}\rightarrow\mathcal{A}^{sa}$ are introduced as an abbreviation for the expressions
	\begin{multline}\label{freeItoFormula-L0}
		L^0\left[f(X_s)\right] = T_{f^{[1]}}^{X_0,X_0}(a(X_s)) +\\+ \sum_{(i,j)\in \{1,\dots,d\}^2}T_{f^{[2]}}^{X_0,X_0,X_0}(b^i(X_s)dW_sc^i(X_s),b^j(X_s)dW_sc^j(X_s))
	\end{multline}
	and
	\begin{equation}\label{freeItoFormula-L1}
		L^1[f(X_s)]_0^t = \sum_{i=1}^d  \int_0^t  T_{f^{[1]}}^{X_0,X_0}(b^i(X_s)dW_sc^i(X_s)).
	\end{equation}
\end{thm}
For the proof we refer to \cite{Biane1998}. In \cite{Nikitop-Ito}, the theorem was derived directly from a Taylor expansion of operator function $f$.\\
In the context of stochastic differential equations, we will make use of  \Cref{freeItoTheorem}  to work out an iterated It\^{o}-Taylor expansion of the solution of the fSDE (see introduction of \Cref{sec:fEMMinLp} and  \Cref{sec:defoffSMandThm}).\\
The It\^{o} formula can also be expressed in product form, see \cite[Theorem 4.1.2]{Biane1998}, \cite{kargin}:
\begin{equation}\label{eq:ItoProductIntegral}
	\begin{aligned}
		\int_0^1a_tdW_tb_t\int_0^1c_tdW_tdt
		&=\int\left( \int_0^t a_sdW_sb_s\right)dW_tc_tdW_td_t\\
		&+\int a_tdW_tb_T \left( \int_0^t a_sdW_sb_s\right)\\
		&+\int_0^1 \varphi(b_tc_t)a_td_tdt.
	\end{aligned}
\end{equation}
This formula will play a major role in developing first order methods. The rule can also be written in differential form as (see \cite{kargin})
\begin{equation}\label{freeItoFormulaAsInKargin}
	a_tdW_tb_t \cdot c_tdW_td_t= \varphi\left(b_tc_t\right)a_td_tdt.
\end{equation}
In the important case $a_t=c_t=d_t=1$ this yields the formal rules 
$$dW_tb_tdW_t = \varphi\left(b_t\right)dt,\quad dW_tdW_t=dt.$$
The formula also implies the formal rules $$dt^2=dW_tdt=0.$$

\section{Numerical Approximation, Euler Maruyama}\label{sec:fEMMinLp}
Let $0<T<\infty$ and $I=[0,T]$. We assume that the fSDE  \eqref{intro-freeSDE-diffform} has a unique solution $(X_t)_{t\in[0,T]}$ adapted to the filtration $\mathbb{F}$. For $L\in\N$, consider a discretization of $I$ as $t_0=0<t_1<\cdots<t_L=T$ and $\Delta t=T/L$.
We recurrently calculate a numerical approximations $(\oX_k)_{k=0,\dots,L}$,  $\oX_k\in\mathcal{A}_{t_k}$,  by 
$$
\oX_{k+1}=\oX_k + \Phi\left(t_k,\oX_k, \Delta t, W_{t_{k+1}}-W_{t_k}\right),
$$
starting with $\oX_0=X_0\in\mathcal{A}^{sa}$. $\Phi:\R\times \mathcal{A} \times \R\times \mathcal{A} \rightarrow \mathcal{A}$ is the increment function, which determines the numerical algorithm. 
In the following we simplify the notation, and write $X_k$ as the solution of \eqref{intro-freeSDE-diffform} evaluated at time point $t_k$ i.e. $X_k=X_{t_k}$. Analog for $W_{t_k}=W_k$, $\Delta W_k=W_{{k+1}}-W_{k}$.

\begin{defn}\label{def-strong-converg-fEMM}
	The numerical approximation $(\oX_k)_{k=0,\dots,L}$ is said to converge strongly to the solution $(X_t)_{t\in I}$ of \eqref{intro-freeSDE-diffform} in $L_p(\varphi)$-norm ($1\le p\le \infty$) with order $\gamma>0$, if there is a constant $C>0$ independent of $\Delta t$, such that
	\begin{equation}
		\|\oX_{k}-X_k\|_p\leq C (\Delta t)^\gamma
	\end{equation}
	for all $k=0,\dots,L$ and $L\in\N$. 
\end{defn}
In  \cite{SchlueWib2023}, the authors introduced the free analog of the Euler-Maruyama method, and proved strong convergence in $L_2(\varphi)$ of order $\gamma=1/2$ (together with weak convergence of order $\gamma=1$). In this paper, we extend this result to strong convergence in  $L_\infty(\varphi)$ with order $\gamma=1/2$. This result will follow immediately from considerations towards methods of order $\gamma=1$.  Nevertheless, we repeat the definition of the free Euler-Maruyama (fEMM) in the following for the sake of completeness. In \Cref{sec:Numeric:Taylorfirst}, we perform the first step of the stochastic It\^{o}-Taylor expansion, already carried out in \cite{SchlueWib2023}, but adding \Cref{lem:estim-ito-taylor-expansion}. In \Cref{sec:Numeric:strongfEMM} the method fEMM will be defined together with the statement of the convergence result in $L_\infty(\varphi)$.
\subsection{Stochastic It\^{o}-Taylor Expansion - First Step}\label{sec:Numeric:Taylorfirst}

Assuming $a,b^i,c^i\in W_3(\R)$ in \eqref{intro-freeSDE-diffform}, we can apply the free It\^{o} formula \eqref{eq:freeItoFormula} for $f=a,b^i,c^i$. This yields an iterated free It\^{o} formula which allows to motivate and define a free analog of the Euler-Maruyama method (see also \cite{SchlueWib2023}). Using the abbreviations $a(X_t)=a_t$ (similar notation for $b^i,c^i$), we obtain
\begin{multline}\label{eq:h1231}
	X_{t+\Delta t}-X_t=\int_t^{t+\Delta t} a_t ds+
	\int_t^{t+\Delta t}\int_t^s L^0[a_u]du\,ds+\int_t^{t+\Delta t}L^1[a_u]_t^{s}ds + \\
	\sum\limits_{i=1}^d\int\limits_t^{t+\Delta t}\left( b_t^i+\int_t^sL^0[b_u^i]du+L^1[b_u^i]_t^{s}) \right)dW_s\left( c_t^i+\int_t^sL^0[c_u^i]du+L^1[c_u^i]_t^{s}) \right)
\end{multline}
Since $a_t,b_t^i,c_t^i$ do not depend on the integration variable $s$, we rewrite \eqref{eq:h1231} as
\begin{equation}\label{eq:h4}
	X_{t+\Delta t}-X_t=a_t\Delta t+\sum_{i=1}^db_t^i(W_{t+\Delta t}-W_t)c_t^i + \sum_{i=1}^d M^i_t(\Delta t) + \sum_{i=1}^dR_t^i(\Delta t),
\end{equation}
where
\begin{equation}\label{eq:h5-1}
	M_t^i(\Delta t) = \int_t^{t+\Delta t} b_t^i dW_s\left(L^1[c_u^i]_t^{s} \right)+
	\int_t^{t+\Delta t}\left(L^1[b_u^i]_t^{s}\right)dW_sc_t^i 
\end{equation}
and
\begin{multline}\label{eq:h5}
	R_t^i(\Delta t) = \int_t^{t+\Delta t}\int_t^s L^0[a_u^i]du\,ds+ \int_t^{t+\Delta t}L^1[a_u^i]_t^sds + \\
	+\int_t^{t+\Delta t} b_t^i dW_s\left( \int_t^sL^0[c_u^i]du\right) + 
	\int_t^{t+\Delta t}\left(\int_t^sL^0[b_u^i]du\right) dW_s c_t^i + \\
	+\int_t^{t+\Delta t}\left(\int_t^sL^0[b_u^i]du\right) dW_s\left( \int_t^sL^0[c_u^i]du\right)+\\
	+\int\limits_t^{t+\Delta t}\left(\int_t^sL^0[b_u^i]du\right) dW_s\left( L^1[c_u^i]_t^{s} \right)
	+\int\limits_t^{t+\Delta t}\left(L^1[b_u^i]_t^{s}\right)dW_s\left( \int_t^sL^0[c_u^i]du\right) + \\
	+\int_t^{t+\Delta t}\left(L^1[b_u^i]_t^{s}\right)dW_s\left( L^1[c_u^i]_t^{s} \right).
\end{multline}

Due to the boundedness and continuity  of the involved functions $a,b^i,c^i$ the above integrals are well defined.
\begin{lem}\label{lem:estim-ito-taylor-expansion}
 Assume $a,b^i,c^i$ are locally operator Lipschitz and $\|X_t\|\leq M$ for $t\in[0,T]$. Then the following estimations hold:
	\begin{align*}
		\left\|\int_t^{t+\Delta t}a(X_s)ds\right\|&=\mathcal{O}(\Delta t)\\
		\left\|\int_t^{t+\Delta t} b(X_s)dW_sc(X_s)\right\|&=\mathcal{O}(\sqrt{\Delta t})\\
		\|M_t^i(\Delta t)\|^2&=\mathcal{O}(\Delta t^2)\\
		\|R_t^i(\Delta t)\|^2&=\mathcal{O}(\Delta t^3).
	\end{align*}
\end{lem}
\begin{proof}
	The first estimation is clear. The second estimation is a direct consequence of the free Burkholder-Gundy inequality (\cite{Biane1998}). The third estimation follows from the iterated stochastic integral in $M_t^i$, Burkholder-Gundy, Lipschitz property of $b,c$, \cite[Theorem 3.5.3]{Skripka2019} and the uniform bound of $\|X_t\|$. Since $\|\Delta W\|^2=\mathcal{O}(\Delta t)$, the iterated integrals are of order $\Delta t$ in operator norm.
	The "smallest" summand in $R_t^i$ is the integral $\int_t^{t+\Delta t}L_1[a_u^i]_t^sds$, which is an iterated integral of a Bochner integral and a free stochastic integral. It is of order $\mathcal{O}(\Delta t^{3/2})$. 
\end{proof}

\subsection{Strong Convergence of the free Euler-Maruyama Method}\label{sec:Numeric:strongfEMM}
The free Euler-Maruyama method can now be motivated from \eqref{eq:h4} simply by  skipping the terms $M_t^i$ and $R_t^i$. The following is taken from \cite{SchlueWib2023}. 
\begin{defn}[fEMM]\label{def:freeEM-Definition}
	Given $T>0$, consider a partition of $[0,T]$ into $L\in\N$ intervals $[t_{k},t_{k+1}],k=0,\dots,L-1$ with constant step size $\Delta t=\frac{T}{L}$.
	Define the one-step free Euler-Maruyama approximation (fEMM) $\oX_k$ of the solution $X_{t_k}=X_k$ of \eqref{intro-freeSDE-diffform} at time point $t_k$ by
	\begin{equation}\label{kap2-def-fEMM}
		\oX_{k+1}=\oX_{k}+a(\oX_{k})\Delta t+\sum_{i=1}^db^i (\oX_{k})\Delta W_{k}c^i(\oX_{k}),\,\,\, 	k=0,1,\dots,L-1
	\end{equation}
	with start value $\oX_0=X_0\in\mathcal{A}^{sa}$. $\oX_k$ denotes the numerical approximation of $X_t$ at time point $t=t_k$. $\Delta W_k$ is the increment of the Brownian motion in \eqref{eq:h4} evaluated at the discretization time points, i.e. $\Delta W_{k}=W_{k+1}-W_{k}=W_{t_{k+1}}-W_{t_k}$.
\end{defn}
In \cite{SchlueWib2023} the authors proved strong convergence order of $\gamma=\frac{1}{2}$ of fEMM \eqref{kap2-def-fEMM} in $L_2$-norm under certain assumptions. In \cite{NiuWeiYinZhen} the authors developed implicit methods supplemented by stability analysis.\\
We are now ready to formulate the following theorem, which extends the $L_2(\varphi)$ convergence of fEMM to $L_\infty(\varphi)$ in the strong sense. For better readability, we formulate the theorem already at this point, although it is a direct consequence of \Cref{thm:thm-order1-general}. We refer the reader to \Cref{sec:MilsteinTypeOrder1} for details.
\begin{thm}\label{thm:fEMMstronginLinfty}
	The free Euler-Maruyama method fEMM \eqref{kap2-def-fEMM} is convergent in the strong sense in $L_p,\, 1\leq p \leq \infty$ with order $\gamma=1/2$. 
\end{thm}
\begin{proof}
The method fEMM is obtained by simply skipping terms $M_t^i$ and $R_t^i$ in \eqref{eq:h4}, i.e. they are not discretized. Then, the statement is a direct consequence of \Cref{thm:thm-order1-general}.
\end{proof}
We refer to \Cref{cor:thm1-FEMM-order-1} which states strong convergence of fEMM of order $\gamma=1$ in case the diffusion is constant, i.e. for equations of the form  $dX_t=a(X_t)dt+ \lambda dW_t, \lambda \in \R $.

\section{Approximations of Strong Order \texorpdfstring{$\gamma=1$ \\ in \(L_p(\varphi), 1\leq p \leq \infty\)}{Approximation of Strong Order gamma}}\label{sec:MilsteinTypeOrder1} 
Starting point is the expansion \eqref{eq:h4}, which includes \eqref{eq:h5} and \eqref{eq:h5-1}. By iterating the expansion (in \Cref{sec:itotaylor-seconditeration}), we will obtain those terms, which must be discretized in order to obtain a method with $\gamma=1$. We will see, that these terms contain iterated free stochastic integrals, just as in the commutative case. A proper discretization of the iterated integrals is required, which will be content of \Cref{sec:num-approx-m-terms}. Unfortunately, it turns out, that only in very special cases, the iterated integrals can be calculated directly. Only the cases
\begin{itemize}
	\item single affine diffusion term ($d=1$) in general $\mathcal{A}$,
	\item single nonlinear diffusion term in the special von Neumann algebra $\mathcal{A}=\mathcal{M}^N(\R)_{sa}$,
\end{itemize}
end up in a derivative free Milstein-type method for fSDEs. All the other cases do rely on a approximation of the iterated integrals. In the case of multiple nonlinear diffusion terms, the iterated integrals in $m_t(\Delta t)$ cannot be resolved by It\^{o} (see \Cref{sec:constmethods-nonlindiff-subdiv}). Up to now, the general case $d>1$ can only be handled by approximating the iterated integrals via the method of subdividing the discretization interval.
\subsection{Stochastic It\^{o}-Taylor Expansion - Second Iteration}
\label{sec:itotaylor-seconditeration}
As in the commutative  case, the key in the development of higher order methods is the iterated stochastic It\^{o}-Taylor expansion of the solution of the underlying differential equation and its discretization. To obtain such an expansion, one needs to apply the free It\^{o} formula \eqref{eq:freeItoFormula} repeatedly to the drift and diffusion terms $a, b^i, c^i$ of the underlying fSDE. In the last chapter we defined fEMM out of the expansion \eqref{eq:h4}.
To be able to define a method of strong order $\gamma=1$, we apply  \eqref{eq:freeItoFormula} to $a,b^i,c^i$ in \eqref{eq:h5-1} once more, but we only take the terms $b_t^i$ and $c_t^i$, which are the functions $b^i,c^i$ evaluated at time point $t$, for example $b^i(X_t)=b_t^i$. This yields
\begin{equation}\label{eq:proof-strong-conv-start-def}
	X_{t+\Delta t}=X_t + a_t\Delta t + \sum\limits_{i=1}^db_t^i\Delta W_t c_t^i + \sum\limits_{i=1}^dm_t^i(\Delta t)+ \sum\limits_{i=1}^d\rho_t^i(\Delta t),
\end{equation}
where the term $m_t^i(\Delta t)$ contains the iterated free stochastic integrals, i.e.
\begin{multline}\label{E:m1rep1-1-def}
	m_t^i(\Delta t) = \int_{t}^{t+\Delta t} b_t^i dW_s\left( 	T_{c^{i,[1]}}^{X_t,X_t}\left(\sum\limits_{j=1}^{d}b_t^j \int_{t}^s dW_uc_t^j \right)\right) + \\ + 
	\int_{t}^{t+\Delta t}\left(T_{b^{i,[1]}}^{X_t,X_t}\left(\sum\limits_{j=1}^{d}b_t^j\int_{t}^sdW_uc_t^j\right)\right)dW_s c_t^i.
\end{multline}
By the same arguments as in  \Cref{lem:estim-ito-taylor-expansion}, we have $\|m_t^i(\Delta t)\|^2=\mathcal{O}(\Delta t^2)$. Note that $m_t^i(\Delta t)$ are the only terms of $\mathcal{O}(\Delta t^2)$ in the expansion \eqref{eq:proof-strong-conv-start-def}. All higher order terms (including $R^i$ in \eqref{eq:h4}) are collected in 
 $\rho_t^i(\Delta t)$. It follows  $\|\rho_t^i(\Delta t)\|^2 = \mathcal{O}(\Delta t^3)$ by the same arguments as for $m_t^i$ above.

\subsection{General Theorem of Strong Convergence with Order \texorpdfstring{$\gamma=1$ in $L_p, \,1\leq p\leq \infty$}{General Theorem of Strong convergence in Lp}}

Let  $t_0=0<t_1<\dots<t_L=T$ a discretization of the interval $[0,T]$. A numerical method based on the It\^{o}-Taylor expansion is build up by finding a proper discretization of the terms in \eqref{E:m1rep1-1-def}, denoted by $\om_k^i(\Delta t)$. The construction of such approximations is subject of \Cref{sec:num-approx-m-terms}. A numerical approximation $\oX_k$ to $X_{t_k}$ is then calculated by
 \begin{equation}\label{E:Milstein-Method-General}
 	\oX_{k+1}=\oX_k+a(\oX_k)\Delta t + \sum_{i=1}^d b^i(\oX_k)\Delta W_kc^i(\oX_k) + \sum_{i=1}^{d}\om_k^i(\Delta t).
 \end{equation}
We will now formulate a general theorem, which imposes conditions on $\om_k^i(\Delta t)$ in order to obtain a convergent method in the strong sense with order $\gamma=1$.
\begin{thm}\label{thm:thm-order1-general}
	Let $a:\mathcal{A}\rightarrow \mathcal{A}$ be an operator function with $a\in W_3(\R)$. The functions $b^i, c^i$ have analog properties as $a$. Assume that there is a constant $\overline{M}>0$ such that $\|\oX_k\|<\overline{M}$ for $k=0,\dots,L$ and $L\in\N$, i.e.  $\overline{M}$ is independent of the discretization.
	Given a numerical method \eqref{E:Milstein-Method-General} for which
	\begin{enumerate}
		\item $\om_k^i(\Delta t)$ belong to the same subalgebra as $m_k^i(\Delta t)$, i.e. \eqref{E:m1rep1-1-def} for at $t=t_k$,

		\item $\label{E:Estim-of-m-order-1}
			\|m_k^i(\Delta t)-\om_k^i(\Delta t)\|^2 \leq C_m\Delta t\|X_k-\oX_k\|^2 + \mathcal{O}(\Delta t^3)$,
	\end{enumerate}
then the numerical approximation converges in the strong sense with order $\gamma=1$, i.e. there is a constant $C_m>0$ such that
	\begin{equation}
		\|X_k-\oX_k\|_p\leq C_m \Delta t
	\end{equation}
	for all $L\in\N$, $k=0,\dots, L$ and $1\leq p\leq \infty$. 
	The constant $C_m>0$ is independent of step size $\Delta t=T/L$ (resp. $k$).
\end{thm}

\begin{proof} Since $a$ is continuous, $a(\mathcal{A}_{sa})\subseteq \mathcal{A}_{sa}$. Consider a discretization of $[0,T]$ as described above. We first build up a continuous reconstruction of $X_t$ out of the discrete values $\oX_k$ obtained from \eqref{E:Milstein-Method-General}. Let's define the order one reconstruction
	\begin{equation}\label{def:fMM-reconstruct-strong-order-1}
		Z_\tau=Z_{k}+\oa_{k}(\tau-t_k) + \ob_{k}(W_\tau-W_{t_k})\oc_{k} + \sum_{i=1}^d\om_{k}^i(\tau-t_k)
	\end{equation}
	on the interval $[t_{k},t], t_k\leq\tau\leq t_{k+1},~k=0,\dots,n_t-1$. Note that $Z_{t_k}$ is written as $Z_k$ and  $Z_\tau$ coincides with $\oX_k$ at the discretization point $\tau=t_k$, i.e. $Z_k=\oX_k$. 
	Next we analyze the difference
	\begin{multline}\label{est:beweis-milstein-1}
		X_t-Z_t = \sum_{k=0}^{n_t-1} (X_{k+1}-X_k) -  \sum_{k=0}^{n_t-1} (Z_{k+1}-Z_k) + (X_t - X_{n_t}) - (Z_t - Z_{n_t}) = \\ =
		\underbrace{\sum_{k=0}^{n_t-1} \left(a_k-\oa_k\right)\Delta t}_{S_1} +
		\underbrace{\sum_{k=0}^{n_t-1} \sum_{i=1}^d\left( b_k^i\Delta W_k c_k^i - \ob_k^i\Delta W_k \oc_k^i\right)}_{S_2}
		+ \\ +  
		\underbrace{\sum_{k=0}^{n_t-1}\sum_{i=1}^d \left(m_k^i(\Delta t)-\om_k^i(\Delta t)\right)}_{S_3} 
		+  \underbrace{(X_t - X_{n_t})}_{S_4} - \underbrace{(Z_t - Z_{n_t})}_{S_5}+ \underbrace{\sum_{k=0}^{n_t-1}\sum_{i=1}^d\rho_k^i(\Delta t)}_{S_6}
	\end{multline}
	Applying the $\|\cdot\|^2$ to the difference $X_t-Z_t$ followed by the triangle inequality and $(u_1+\dots +u_6)^2\leq 6(u_1^2+\dots u_6^2)$ leaves the task to estimate $\|S_i\|^2,~i=1\dots 6$.
	Now define
	\begin{displaymath}
		v(u)=\sup_{0\leq s\leq u}\left\|X_s-Z_s\right\|^2.
	\end{displaymath}
	Since $a$ is locally operator Lipschitz in $L_\infty(\varphi)$ and $Z_k=\oX_k$ we obtain
	\begin{multline}\label{eq:chap7:estim-term-S1}
		\|S_1\|^2=\left\|\sum_{k=0}^{n_t-1} \left(a_k-\oa_k\right)\Delta t \right\|^2\leq 
		\sum_{k=0}^{n_t-1} n_t L_a \Delta t^2 \|X_k-Z_k\|^2
		\leq \\ \leq
		\sum_{k=0}^{n_t-1} \Delta t K_{11} \|X_k-Z_k\|^2  \leq K_1\int_0^{t_{n_t}}v(u)du.
	\end{multline}
	Note that $\Delta t = t_{k+1}-t_k \geq t-t_k$. In a similar way we estimate the second summand in \eqref{est:beweis-milstein-1}  as
	\begin{equation}\label{eq:chap7:estim-term-S2}
		\|S_2\|^2	\leq K_{2}\int_0^{t_{n_t}}v(u)du.
	\end{equation}
	Since $\|\rho_k^i(\Delta t)\|^2=\mathcal{O}(\Delta t^3)$ we get
	\begin{equation}\label{eq:chap7:estim-term-S6}
		\|S_6\|^2	\leq K_{6}\Delta t^2.
	\end{equation}
	Now we turn to the estimation of the terms $S_3$ and $S_6$. We strongly depend on a result of D.~Voiculescu (\cite{VoicuAddRV}). If $X_i\in\mathcal{A}, i=1,\dots n$ are free random variables with $\varphi(X_i)=0$, then $\|\sum_{i=1}^n X_i \|\leq \sup_i \|X_i\| + 2\sqrt{\sum_{i=0}^n\|X_i\|_2^2}$. Since $\om_k^i(\Delta t)$ belong to the same subalgebra as $m_k^i(\Delta t)$, the differences $m_k^i(\Delta t)-\om_k^i(\Delta t)$ are free with respect to $k$. Then we are able to estimate $S_3$ as follows.
	\begin{align*}\label{eq:chap7:estim-term-S3}
		\|S_3\|^2 &\leq \left\|
		\sum_{i=1}^d\sum_{k=0}^{n_t-1}m_k^i(\Delta t)-\om_k^i(\Delta t)\right\|^2
		\\
		&\leq d\sum_{i=1}^d \left\|\sum_{k=0}^{n_t-1}m_k^i(\Delta t)-\om_k^i(\Delta t)\right\|^2
		\\
		&\leq d\sum_{i=1}^d 
		\left(\sup_{0\leq k\leq n_t-1} \|m_k^i(\Delta t) - \om_k^i(\Delta t)\| + 2\sqrt{\sum_{k=0}^{n_t-1}\|m_k^i(\Delta t)-\om_k^i(\Delta t)\|_2^2}\right)^2\\
		&\leq d\sum_{i=1}^d 
		\left(\sup_{0\leq k\leq n_t-1} \|m_k^i(\Delta t) - \om_k^i(\Delta t)\| + 2\sqrt{\sum_{k=0}^{n_t-1}\|m_k^i(\Delta t)-\om_k^i(\Delta t)\|^2}\right)^2\\
	\end{align*}
	Since $(\sup_i a_i)^2=\sup_i a_i^2$ for non negative $a_i, i=0,\dots,n$  
	and assumption \eqref{E:Estim-of-m-order-1} we proceed as follows.
	\begin{align*}
		\|S_3\|^2 &\leq d\sum_{i=1}^d \left(
		\sqrt{\sup_{0\leq k\leq n_t-1} \|m_k^i(\Delta t) - \om_k^i(\Delta t)\|^2} +  2\sqrt{\sum_{k=0}^{n_t-1}\|m_k^i(\Delta t)-\om_k^i(\Delta t)\|^2}\right)^2\\
		&\leq d\sum_{i=1}^d \left(3
		 \sqrt{\sum_{k=0}^{n_t-1}\|m_k^i(\Delta t)-\om_k^i(\Delta t)\|^2}\right)^2\\
		&= 9d\sum_{i=1}^d 
		\sum_{k=0}^{n_t-1} \|m_k^i(\Delta t) - \om_k^i(\Delta t)\|^2\\
		&\leq 9d\sum_{i=1}^d 
		\sum_{k=0}^{n_t-1} C_m\Delta t\|X_k-\oX_k\|^2 + 		\mathcal{O}(\Delta t^3)\\
		&\leq K_1 \int_0^{t_{n_t}}v(u)du + K_2\Delta t^2.
	\end{align*}
Similar arguments as in the estimation of $S_1, S_2, S_3$ and $S_6$ lead to 
	\begin{equation}\label{eq:chap7:estim-term-S45}
		\|S_4-S_5\|^2\leq K_{45}\int_{n_t}^tv(u)du + C_{45}\Delta t^3.
	\end{equation}
Collecting \eqref{eq:chap7:estim-term-S1} to \eqref{eq:chap7:estim-term-S45} allows to estimate 
	\begin{equation}\label{eq:proof-strong-fMM-before-gronwall}
		v(t)\leq C_1\int_0^tv(u)du + C_2\Delta t^2.
	\end{equation}
	A Gronwall argument results in the following estimation of $v$,
	\begin{equation*}
		v(t)=\sup_{0\leq s\leq t}\|X_s-Z_s\|^2\leq C_3\Delta t^2.
	\end{equation*}
	This implies
	\begin{equation}
		\| X_k-\oX_k\|\leq C\Delta t,
	\end{equation}
	for all $k=0,\dots,L$. Since the trace is unital it follows
	\begin{equation}
		\| X_k-\oX_k\|_p\leq C\Delta t
	\end{equation}
	for $1\leq p \leq \infty$.
\end{proof}

In the special case of constant diffusion the term \eqref{E:m1rep1-1-def} vanishes and it follows
\begin{cor}\label[corollary]{cor:thm1-FEMM-order-1}
	Consider the fSDE
	$$
	dX_t=a(X_t)dt+ \lambda dW_t.
	$$
	Then the  free Euler-Maruyama method fEMM \eqref{kap2-def-fEMM} is strong convergent in $L_p,\, 1\leq p \leq \infty$ with order $\gamma=1$. This extends the result of \cite{SchlueWib2023}.
\end{cor}
\subsubsection{Boundedness of the Numerical Solution}
In \Cref{thm:thm-order1-general} we assumed a uniform bound on the numerical solution. We now show, that if a numerical method  performs with strong convergence of order $\gamma>0$, this assumption is true.
\begin{thm}\label{thm:strong-conv-bound-num-sol}
	Consider a fSDE \eqref{intro-freeSDE-diffform} with the solution $X_t$ on $[0,T]$ where $T<\infty$ (see \Cref{rem:fSDE-sol-theory}).
	Given a discretization $T=L\Delta t$ as defined in \Cref{def:freeEM-Definition}, the numerical solution $\oX_k$ is calculated via a method of first order as defined in \eqref{E:Milstein-Method-General}. Let $a$, $b$ and $c$ be operator functions which are locally operator Lipschitz as in \Cref{thm:thm-order1-general}.  Then the numerical solution is uniformly bounded for each $k=0,\dots,L, \, L\in\N$, i.e. there is a constant $\tilde{M}>0$ such that $\|\oX_k\|<\tilde{M}$, where $\tilde{M}$ does not depend on $L$ resp. $\Delta t$.\\
\end{thm}
\begin{proof}
	For the following, we construct a piecewise constant process $(\oX_t)_{t\in[0,T]}$ defined by $\oX_t=\oX_k$ for $t_k\leq t<t_{k+1}$, where $\oX_k$ is a numerical solution on $[0,T]$ calculated by a step size $\Delta t$. 
	We have $\|X_0\|<M$. Let $\epsilon>0$ independent of the discretization. Therefore, for each $\Delta t=T/L$, there is a time point $0< T^{(\Delta t)}\leq T$, such that $\|X_t-\oX_t\|\leq C\Delta t$ for all $t\in [0,T^{(\Delta t)}]$ and $T^{(\Delta t)}=\sup\{0<\tilde{T}\leq T, \, \sup_{t\in[0,\tilde{T}]}\|\oX_t\|<M+\eps\}$. From this it follows that  $\|\oX_{t}\|\leq M+C\Delta t$ for all $t\in[0,T^{(\Delta t)}]$.
	Let $l(\Delta t)\in\N$ such that $t_{l(\Delta t)}\leq T^{(\Delta t)}<t_{l(\Delta t)+1}$.	
	Applying the triangle inequality to \eqref{E:Milstein-Method-General} yields (using the abbreviation $l=l(\Delta t)$)
	\begin{equation*}
		\|\oX_{l+1}\| \leq  \|\oX_{l}\|+ \|a(\oX_{l})\|\Delta t + \sum_{i=1}^d \left(\|b^i(\oX_{l})\Delta W_lc^i(\oX_l)\| + \|\om_{l}^i(\Delta t)\|\right)
	\end{equation*}
	Since $\|\oX_{l}\|\leq M+C\Delta t$, the Lipschitz property of $a$ and Lemma \ref{lem:discstochint_lin_growth} allow to estimate 
	\begin{equation*}
		\|\oX_{l+1}\| \leq  M+\mathcal{O}(\sqrt{\Delta t}) + \|\om_l^i(\Delta t)\|.
	\end{equation*}
	Due to the second assumption in \Cref{thm:thm-order1-general}   and the bounds on $X_l$ and $\oX_l$, we have $\|\om_l^i(\Delta t)\|=\sqrt{\mathcal{O}(\Delta t)}$.
	If $\Delta t$ is small enough, it follows
	$$
	\|\oX_{l+1}\|< M+\epsilon.
	$$
	But this is a contradiction to  $T^{(\Delta t)}$ as the supremum over all time points in $[0,T]$, for which the numerical solution is bounded by $M+\epsilon$. This shows, that for each $\epsilon>0$ and $\Delta t$ small enough, $T^{(\Delta t)}=T$. Therefore there is a $\tilde{M}>M$ such the numerical solution is bounded by $\tilde{M}$ in operator norm, independent of the discretization.
\end{proof}

\section{Construction of Methods}\label{sec:num-approx-m-terms}
The development of a Milstein type method of order $\gamma=1$ requires the discretization of the iterated free stochastic integrals \eqref{E:m1rep1-1-def}. It is natural to ask, under which circumstances these integrals can be calculated directly. As it turns out, there's a major difference between $d=1$ and $d>1$. In the first case, It\^{o} can be applied, in the second case not. Therefore this section starts with a discussion of the case $d=1$. Due to the imposed non-commutativity the treatment of  \eqref{E:m1rep1-1-def} differs from the commutative case even for a single diffusion $d=1$. We have to deal with additional triple operator integrals, which stand in place for the second derivative of the diffusion operator functions. Unfortunately theses terms are only of order $\Delta t^2$ (in square of $L_\infty(\varphi))$-norm, which does not allow to skip them in the development of first order methods. These terms do not appear in the commutative case, where the classical Milstein method only contains first order derivatives. \\
In \Cref{sec:dequal1} we discuss, under which conditions these operator integrals can be calculated directly in case of a single diffusion term. Unfortunately, we will see in section \Cref{subsec:dgt1nonlinsubdiv} that all these methods fail for the case $d>1$. Even more, in \Cref{subsec:dgt1nonlinsubdiv} it is stated, that the It\^{o} formula cannot be applied to resolve the iterated integrals as in the case $d=1$. We therefore present a method to discretize the iterated integrals to achieve order $\gamma=1$. In fact, to the best  of our knowledge, there's no other discretization of iterated integrals in the given context known in the literature. 
 \begin{rem}
 	In \cite{ortmann2011functionalsfreebrownianbridge}, introduced the notion of a Free Brownian Bridge, which could possibly be an attempt to approximated the iterated free stochastic integrals.
 \end{rem}
\subsection{\texorpdfstring{Methods for Single Diffusion, $d=1$}{Methods for Single Diffusion, d}}\label{sec:dequal1}
Now we turn to the construction of methods \eqref{E:Milstein-Method-General} in the simplest case of a single diffusion term $d=1$. We will see, that due to the non-commutativity, only in special cases a simple discretization can be found, as already mentioned.\\
In the commutative case, the key to the construction of higher order numerical approximation relies on the discretization of the iterated stochastic integrals in the stochastic Taylor expansion of the solution of the underlying SDE (e.g. see \cite{kloedenplaten}, \cite{milsteintretyakov}). Consider the simplest case of a scalar SDE with a single Brownian motion and differentiable diffusion term $b$, i.e.
$$
dx_t=a(t,x_t)dt + b(t,x_t)dW_t.
$$
The stochastic Taylor expansion for a method of strong order $1$ includes the following iterated integral, which can be integrated exactly, i.e.
$$
\int_{t_k}^{t_{k+1}}\int_{t_k}^s dW_u dW_s =\frac{1}{2}\left((\Delta W_{k})^2-\Delta t\right).
$$
The resulting method is the well-known Milstein method (\cite{milsteintretyakov}, \cite{kloedenplaten}). For SDE with more complicated structure, e.g. multiple Brownian motions, an efficient evaluation of the iterated integrals is only possible, if the structure of the SDE has special properties, i.e. commutative noise.\\
In the non-commutative case, the structure of the terms in the stochastic Taylor expansion \eqref{E:m1rep1-1-def} consists of a sum of iterated integrals. The goal is, to combine these summands in order to convert them into a product by applying the It\^{o} formula in product form, \cite[Theorem 4.1.2]{Biane1998}.\\
It can easily be seen in \eqref{E:m1rep1-1-def} that the application of the product form of the It\^{o}  formula requires $c^i=b^j, i,j=1,\dots d$. This means $d=1$, which implies $b=c$. In this case \eqref{E:m1rep1-1-def} simplifies to
\begin{multline}\label{E:fMMvariante1-def-1}
	m_t(\Delta t)= 
	b_t\left(\int_{t}^{t+\Delta t} dW_sT_{b^{[1]}}^{X_t,X_t}(  \int_{t}^s b_tdW_u) \right.
	+\\+
	\left.
	\int_{t}^{t+\Delta t} T_{b^{[1]}}^{X_t,X_t}(\int_{t}^s dW_u b_t )dW_s
	\right)b_t.
\end{multline}
Note that we applied the fact, that two elements $U,V\in\mathcal{A}$ commute, if both belong to the same subalgebra, e.g. $U,V\in \langle X_t\rangle \subset \mathcal{A}$.\\
In order to utilize the It\^{o} formula we need to commute $dW_s$ into the operator integrals, such that the outer stochastic integrals are pushed into the argument of the double operator integral $T_{b^{[1]}}^{X_t,X_t}(\cdot)$. In the non-commutative case, this can be achieved by \Cref{lem:referee}, but with extra costs. Applying  \Cref{lem:referee} to \eqref{E:fMMvariante1-def-1} yields
\begin{multline}\label{E:fMMvariante1-referee}
	m_t(\Delta t)= 
	b_t\left(T_{b^{[1]}}^{X_t,X_t}(  \int_{t}^{t+\Delta t} dW_s\int_{t}^s b_tdW_u) + \right.
	\\+
	\left.
	T_{b^{[1]}}^{X_t,X_t}(\int_{t}^{t+\Delta t}\int_{t}^s dW_u b_tdW_s )
	\right)b_t
	+\\+
	b_t\left(T_{b^{[2]}}^{X_t,X_t,X_t}(V_t,U_t^l) - 
	T_{b^{[2]}}^{X_t,X_t,X_t}(U_t^r,V_t)\right)b_t,
\end{multline}
where we shortened the notation by $V_t=\int_t^{t+\Delta t}dW_sX_t-X_t\int_t^{t+\Delta t}dW_s$,  and $U_t^l=\int_t^{t+\Delta t} b_t dW_u$, $U_t^r=\int_t^{t+\Delta t}dW_ub_t$. The first two summands in \eqref{E:fMMvariante1-referee} are ready to be handled by It\^{o}, i.e.
\begin{multline}\label{E:fMMvariante1-def-2x}
	m_t(\Delta t)=
	b_t\left(T_{b^{[1]}}^{X_t,X_t}\left( \int_{t}^{t+\Delta t} dW_s b_t \int_t^{t+\Delta t}dW_s-\int_t^{t+\Delta t}\varphi(b_t)\textbf{1}dt\right)\right)b_t
	+\\+
	b_t\left(T_{b^{[2]}}^{X_t,X_t,X_t}(V_t,U_t^l) - 
	T_{b^{[2]}}^{X_t,X_t,X_t}(U_t^r,V_t)\right)b_t
\end{multline}
Due to the additional triple operator integrals $T_{b^{[2]}}^{X_t,X_t,X_t}$ we face new difficulties. Unfortunately due to \Cref{lem:estim-ito-taylor-expansion}, the two extra terms obey the property
$$
\|T_{b^{[2]}}^{X_t,X_t,X_t}(V_t,U_t^l)\|^2=\mathcal{O}(\Delta t^2),\,\|T_{b^{[2]}}^{X_t,X_t,X_t}(U_t^r,V_t)\|^2=\mathcal{O}(\Delta t^2).
$$
This shows, that these terms cannot be dropped. They request for a discretization in order to achieve a method of strong order $\gamma=1$.\\

\subsubsection{\texorpdfstring{Linear Diffusion $b$, $d=1$}{Linear Diffusion b}}
Assume that the diffusion in the fSDE is linear, i.e. $b(X_t)= X_t$ (a factor $\mu\neq 0$ will be imposed in the fSDE directly, see \Cref{def:fMMlinearDiff}). This allows a direct calculation of the operator integrals in \eqref{E:fMMvariante1-def-2x}. If $b$ is affine, then $T_{b^{[1]}}^{X,X}(Y)=Y , \, X,Y\in\mathcal{A}$. Since the divided difference $b^{[2]}=0$, the terms including $T_{b^{[2]}}^{X_t,X_t,X_t}$  in \eqref{E:fMMvariante1-def-2x} vanish. We then obtain
\begin{equation*}
	m_t(\Delta t)=X_t\left(\Delta W_t X_t \Delta W_t -  \varphi(X_t)\Delta t\right) X_t.
\end{equation*}
The discretization is the simply constructed via
\begin{equation}\label{eq:deq1lindiffm}
	\om_k(\Delta t)= \oX_k\left(\Delta W_k \oX_k \Delta W_k -  \varphi(\oX_k)\Delta t\right) \oX_t.
\end{equation}
In detail:
\begin{thm}\label{def:fMMlinearDiff}
	Given $T>0$ and the fSDE \eqref{intro-freeSDE-diffform} with $d=1$ and $b_t=c_t$ and $b(X_t)=X_t$.  Consider the fSDE
	\begin{equation*}
		dX_t=a(X_t)dt + \mu X_t dW_t X_t, t\in[0,T],
	\end{equation*}
	with start value $X_0\in\mathcal{A}^{sa}$, $\mu\neq 0$ and a discretization of $[0,T]$ similar to \Cref{thm:strong-conv-bound-num-sol}.
	Then the one-step free Milstein approximation $\oX_k$ to the solution $X_{t}$ on $[0,T]$ is defined via \eqref{eq:deq1lindiffm}, through
	\begin{equation}\label{eq:def-method-linear-d1}
		\oX_{k+1}=\oX_{k}+a(\oX_{k})\Delta t+ \mu \oX_{k}\Delta W_{{k}}\oX_{k}
		+\om_k(\Delta t).
	\end{equation}
	The method is strongly convergent with order $\gamma=1$.
\end{thm}
\begin{proof} Let $Y_k=\Delta W_kX_k\Delta W_k - \varphi(X_k)\Delta t$ and $\oY_k=\Delta W_k\oX_k\Delta W_k-\varphi(\oX_k)\Delta t$. Due to \Cref{thm:thm-order1-general}, we need to estimate 
\begin{align*} 
	m_k(\Delta t) - \om_k(\Delta t)
	&=
	X_kY_kX_k-\oX_k\oY_k\oX_k\\
	&=(X_k-\oX_k)Y_kX_k+\oX_k(Y_kX_k-\oY_k\oX_k)\\
	&=(X_k-\oX_k)Y_kX_k+\oX_k(Y_kX_k-\oY_k X_k)+\oX_k(\oY_k X_k-\oY_k\oX_k)\\
	&=(X_k-\oX_k)Y_kX_k+\oX_k(Y_k-\oY_k)X_k+\oX_k\oY_k (X_k-\oX_k)
\end{align*}
The assumptions state, that $\|X_k\|$ and $\|\oX_k\|$ are uniformly bounded by the constants $M>0$ resp. $\overline{M}>0$. The operator norm is submultiplicative, therefore
\begin{align*}
	\|m_k(\Delta t) - \om_k(\Delta t)\|
	&\leq
	\|X_k-\oX_k\|\|Y_k\|\|X_k\|\\
	&+\|\oX_k\|\|Y_k-\oY_k\|\|X_k\|
	+\|\oX_k\|\|\oY_k\|\|X_k-\oX_k\|\\
	&\leq M\|X_k-\oX_k\|\|Y_k\|+\overline{M}M\|Y_k-\oY_k\| + \overline{M}\|\oY_k\|\|X_k-\oX_k\|.
\end{align*}
Since $\|\Delta W\|=\mathcal{O}(\sqrt{\Delta t})$, $\Delta t<T$ and $\varphi(X)\leq \|X\|_1\leq \|X\|, X\in\mathcal{A}$, we obtain
$\|Y_k\|\leq \|\Delta W^2\|\|X\|+\|X\|\Delta t \leq D_1$ and  $\|\oY\|\leq D_2$.
To estimate $\|Y-\oY\|$, we apply the triangle inequality to get
\begin{align*}
	\|Y_k-\oY_k\|\leq \|\Delta W^2\|\|X_k-\oX_k\|+\varphi(X_k-\oX_k)\Delta t\textbf{1} \leq D_3\|X_k-\oX_k\|.
\end{align*} 
Then 
\begin{equation*}
	\|m_k(\Delta t) - \om_k(\Delta t)\|^2 \leq (3MD_1 + 3\overline{M}MD_3+3\overline{M}D_3)\|X_k-\overline{X}_k\|^2.
\end{equation*}
Since by construction $\om_k(\Delta t)$ are free random variables (with respect to $k$), the statement of the theorem follows from \Cref{thm:thm-order1-general}.
\end{proof}

\subsubsection{\texorpdfstring{Nonlinear Diffusion -  $\mathcal{A}=\mathcal{M}^N(\R)$, $d=1$}{Nonlinear Diffusion - Matrix level}}\label{subsubsec:deq1nonlinmatrix}
For a general von Neumann algebra $\mathcal{A}$ and nonlinear $b$ the discretization of $m_t(\Delta t)$ in \eqref{E:fMMvariante1-def-2x} is unknown, to the best of our knowledge. But if one restricts the considerations to the algebra $\mathcal{A}=\mathcal{M}^N(\R)_{sa}$, then it is possible to calculate the operator integrals by spectral analysis of $X_t$ resp. $\oX_t$ directly. Therefore we can build up a Milstein Method of order $\gamma=1$ based on the stochastic Taylor expansion.\\
In the non commutative probability space $\mathcal{A}=\mathcal{M}^N(\R)_{sa}$, the operator integrals can be expressed as (\cite[Theorem A.9]{Nikitop-Ito}, \cite{azamov_carey_dodds_sukochev_2009}, \cite{Skripka2019}) 
\begin{align}\label{E:Discret-OpInt-Matrix}
	T_{b^{[1]}}^{X,X}(Y)&=\sum_{\boldsymbol{\lambda} \in \sigma(X)^2} b^{[1]}(\boldsymbol{\lambda})P_{\lambda_1}YP_{\lambda_2},\\
	T_{b^{[2]}}^{X,X,X}(Y_1,Y_2) &= \sum_{\boldsymbol{\lambda} \in \sigma(X)^3} b^{[2]}(\boldsymbol{\lambda})P_{\lambda_1}Y_1P_{\lambda_2}Y_2P_{\lambda_3},
\end{align}
where $\sigma(X)$ denotes the spectrum of $x\in\mathcal{A}$ and $P_{\lambda_i}$ are the projection matrices onto the eigenspace of the eigenvalue $\lambda_i\in\R$.
 For readability we write \eqref{E:fMMvariante1-def-2x} slightly shorter.  Evaluated at  $t=t_k$ expression \eqref{E:fMMvariante1-def-2x} reads
\begin{equation}\label{E:fMMvariante1-def-2x-1}
	m_k(\Delta t)=
	b_k\left(T_{b^{[1]}}^{X_k,X_k}\left( Y_k\right)
	+
	T_{b^{[2]}}^{X_k,X_k,X_k}(V_k,U_t^l) - 
	T_{b^{[2]}}^{X_k,X_k,X_k}(U_k^r,V_k)\right)b_k,
\end{equation}
where $Y_k=\int_{t_k}^{t_k+\Delta t} dW_s b_k \int_{t_k}^{t_k+\Delta t}dW_s-\int_{t_k}^{t_k+\Delta t}\varphi(b_{k})ds$. 
The numerical method can be constructed from \eqref{E:fMMvariante1-def-2x-1} via 
\begin{equation}\label{E:m-matrix-method-deq1}
	\om_k(\Delta t)=
	\ob_k \left(
	T_{b^{[1]}}^{\oX_k,\oX_k}(\oY_k)
	+T_{b^{[2]}}^{\oX_k,\oX_k,\oX_k}(\overline{V}_k,\oU_k^l)  - T_{b^{[2]}}^{\oX_k,\oX_k,\oX_k}(\oU_k^r, \overline{V}_k)
	\right)\ob_k,
\end{equation}
where
\begin{align*}
	\ob_k&=b(\oX_{t_k})\\
	\oY_k&=\Delta W_k \ob_k \Delta W_k - \varphi(\ob_k)\Delta t \\
	\overline{V}_k&=\Delta W_k \oX_k - \oX_k \Delta W_k \\
	\overline{U}_k^l &= \ob_k\Delta W_k \\
	\overline{U}_k^r &= \Delta W_k \ob_k \\
	\Delta W_k &= W_{k+1}-W_k.
\end{align*}

\begin{thm}\label{thm:deq1-matrix}
	Let $\mathcal{A}={\mathcal{M}^N(\R)}_{sa}$.  Let the same assumptions as in \Cref{thm:thm-order1-general} be given, except the stronger condition $b\in W_4(\R)$. Take $\om_k(\Delta t)$ from  \eqref{E:m-matrix-method-deq1}. Then the numerical method defined by
	\begin{equation}\label{eq:deq1-matrix}
		\oX_{k+1}=\oX_k+a(\oX_k)\Delta t + b(\oX_k)\Delta W b(\oX_k) + \om_k(\Delta t),
	\end{equation}
	shows strong convergence of order $\gamma=1$.
\end{thm}
\begin{rem} Due to  \Cref{lem:estim-optint-diff-triple}  we have to assume $b\in W_4(\R)$. Since $W_4(\R)\subset W_3(\R)$ the main \Cref{thm:thm-order1-general} can be applied.
\end{rem}
\begin{proof}
	To prove the statement we need to estimate the difference of \eqref{E:fMMvariante1-def-2x-1} and \eqref{E:m-matrix-method-deq1}. We simplify the notation for readability reasons.
	\begin{multline*}
		m_k(\Delta t)-\om_k(\Delta t) = bT_{b^{[1]}}(Y)b-\ob\, T_{b^{[1]}}(\oY)\ob +\\+ b(T_{b^{[2]}}(V,U^l)-T_{b^{[2]}}(U^r,V))b - \ob(T_{b^{[2]}}(\oV,\oU^l) -T_{b^{[2]}}(\oU^r,\oV) )\ob.
	\end{multline*}
	Then 
	\begin{multline}\label{eq:h61}
		\|m_k(\Delta t)-\om_k(\Delta t)\|^2 \leq  2\|bT_{b^{[1]}}(Y)b-\ob\, T_{b^{[1]}}(\oY)\ob\|^2 +\\+ 2\|b(T_{b^{[2]}}(V,U^l)-T_{b^{[2]}}(U^r,V))b - \ob(T_{b^{[2]}}(\oV,\oU^l) -T_{b^{[2]}}(\oU^r,\oV) )\ob\|^2.
	\end{multline}
	We are going to estimate both summands in \eqref{eq:h61} on the right hand side of above inequality. We start with the first one by $bT_{b^{[1]}}(Y)b-\ob\, T_{b^{[1]}}(\oY)\ob = 
	(b-\ob)T_{b^{[1]}}(Y)b + \ob T_{b^{[1]}}(\oY)b-\ob\, T_{b^{[1]}}(\oY)\ob$. Since the norm is submultiplicative (skipping the argument $Y$ and simplify $T_{b^{[1]}}(\oY)=\oT_{b^{[1]}}$)
	\begin{align*}
		\|bT_{b^{[1]}}b-\ob\, \oT_{b^{[1]}}\ob \| &\leq 
		\|b-\ob\|\|T_{b^{[1]}}\|\|b\| + \|\ob\|\| T_{b^{[1]}}b -\oT_{b^{[1]}}\ob \|\\
		&=\|b-\ob\|\|T_{b^{[1]}}\|\|b\| + \|\ob\|\| T_{b^{[1]}}b - T_{b^{[1]}}\ob + T_{b^{[1]}}\ob-\oT_{b^{[1]}}\ob \|\\
		&\leq \|b-\ob\|\|T_{b^{[1]}}\|\|b\| + \|\ob\|\|T_{b^{[1]}}\|\|b - \ob\| + \|T_{b^{[1]}}-\oT_{b^{[1]}}\| \|\ob \|
	\end{align*}
	Due to \cite[Theorem 4.3.8]{Skripka2019}, the local Lipschitz condition on $b$ and the uniform bounds on $\|X_k\|$ and $\|\oX_k\|$ we have
	\begin{equation*}
		|bT_{b^{[1]}}(Y)b-\ob\, \oT_{b^{[1]}}(\oY)\ob \|^2\leq D_1\|X_k-\oX_k\|^2+D_2\|T_{b^{[1]}}-\oT_{b^{[1]}}\|^2.
	\end{equation*}
	To estimate the difference of the operator integrals we first rewrite 
	\begin{multline*}
	T_{b^{[1]}}^{X_k,X_k}\left( Y_k\right) - T_{b^{[1]}}^{\oX_k,\oX_k}(\oY_k)
	=\\=
	T_{b^{[1]}}^{X_k,X_k}\left( Y_k\right) -
	T_{b^{[1]}}^{\oX_k,\oX_k}\left( Y_k\right) +
	T_{b^{[1]}}^{\oX_k,\oX_k}\left( Y_k\right) -
	T_{b^{[1]}}^{\oX_k,\oX_k}(\oY_k).
	\end{multline*}
	By \Cref{lem:estim-optint-diff} it follows that
	\begin{equation*}
				\left\| T_{b^{[1]}}^{X_t,X_t}\left( Y_k\right) - T_{b^{[1]}}^{\oX_k,\oX_k}(\oY_k)\right\|^2
				\leq
				2K_1\|Y\| \|X_k-\oX_k\|^2 + 2\|T_{b^{[1]}}^{\oX_k,\oX_k}(Y_k-\oY_k)\|^2\\
	\end{equation*}
	The double operator integral on the right side can be estimated by \cite[Theorem 4.3.8]{Skripka2019}, assuming a uniform bound on the numerical approximation, i.e.	\begin{multline*}
		\left\| T_{b^{[1]}}^{X_t,X_t}\left( Y_k\right) - T_{b^{[1]}}^{\oX_k,\oX_k}(\oY_k)\right\|^2
		\leq \\ \leq
		K_{11} L_bT^2\|X_k-\oX_k\|^2 + K_{12}\|Y_k-\oY_k\|^2 
		\leq \\ \leq 
		K_{13}\|X_k-\oX_k\|^2 + K_{12}\|X_k-\oX_k\|^2=K_{14}\|X_k-\oX_k\|^2.
	\end{multline*}
	Finally the first summand in \eqref{eq:h61} obeys
	\begin{equation}\label{eq:h62}
		\|bT_{b^{[1]}}(Y)b-\ob\, \oT_{b^{[1]}}(\oY)\ob \|^2\leq D_3\|X_k-\oX_k\|^2.
	\end{equation}
	Now we prepare for the second summand on the right hand side in \eqref{eq:h61}.	The difference of the triple operator integrals in \eqref{E:m-matrix-method-deq1},
	$$
		T_{b^{[2]}}^{X_k,X_k,X_k}(V_k,U_k^l) -
		T_{b^{[2]}}^{\oX_k,\oX_k,\oX_k}(\overline{V}_k,\oU_k^l),
	$$
	 can be estimated by \Cref{lem:estim-optint-diff-triple} as
	 \begin{multline*}
	 \|T_{b^{[2]}}^{X_k,X_k,X_k}(V_k,U_k^l) -
	 T_{b^{[2]}}^{\oX_k,\oX_k,\oX_k}(\overline{V}_k,\oU_k^l)\|
	 \leq \\ \leq
	 K_{21}\|V_k-\oV_k\|+K_{22}\|X_k-\oX_k\|+K_{23}\|U_k^l-\oU_k^l\|,
	 \end{multline*}
	where $V_k-\oV_k=\Delta W_k(X_k-\oX_k)+(X_k-\oX_k)\Delta W_k$ and $U_k^l-\oU_k^l=(b_k-\ob_k)\Delta W_k$. Since $b$ is local operator Lipschitz and $\Delta t<T$ it follows
	\begin{multline*}
		\|T_{b^{[2]}}^{X_k,X_k,X_k}(V_k,U_k^l) -
	T_{b^{[2]}}^{\oX_k,\oX_k,\oX_k}(\overline{V}_k,\oU_k^l)\|^2
	\leq \\ \leq
	(3K_{24}+3K_{25}+3K_{26})\|X_k-\oX_k\|^2 \leq K_{27}\|X_k-\oX_k\|^2.
	\end{multline*}
	Similarly it follows that 
	\begin{align*}
		\left\|T_{b^{[2]}}^{X_t,X_t,X_t}(U_t^r,V_t)-
		T_{b^{[2]}}^{\oX_k,\oX_k,\oX_k}(\oU_k^r, \overline{V}_k)\right\|^2\leq K_{37}\|X_k-\oX_k\|^2.
	\end{align*}
	This allows, with similar estimations performed to obtain \eqref{eq:h61} and \eqref{eq:h62}, that we can finally estimate \eqref{eq:h61} as
	\begin{equation}
		\|m_k(\Delta t)-\om_k(\Delta t)\|^2 \leq (D_3+D_4)\|X_k-\oX_k\|^2.
	\end{equation}
	Then \Cref{thm:thm-order1-general} finishes the proof.
\end{proof}
\begin{rem}
	The previous proof also holds if one drops the restriction on the von Neumann algebra. Since no proper discretization of the operator integrals in a general von Neumann algebra $\mathcal{A}$ is known, we formulated the theorem and proof in $\mathcal{M}^N(\R)_{sa}$.
\end{rem}
\subsubsection{\texorpdfstring{Nonlinear Diffusion -  Subdivision Method, $d=1$}{Nonlinear Diffusion -  Subdivision Method, deq1}}\label{sec:constmethods-nonlindiff-subdiv}
Now we drop the assumption $\mathcal{A}=\mathcal{M}^N(\R)_{sa}$ of \Cref{subsubsec:deq1nonlinmatrix} and consider the fSDE in a general von Neumann algebra $\mathcal{A}$. Here, the problem is that for a method of strong order $\gamma=1$, each of the terms including $T_{b^{[2]}}^{X_t,X_t,X_t}$ in \eqref{E:fMMvariante1-def-2x} have to be discretized. So far, there is no proper discretization of triple operator integrals $T_{b^{[2]}}^{X_t,X_t,X_t}$ in a general von Neumann Algebra $\mathcal{A}$ known (comparable to \Cref{lem:discrete-opt-int}).\\
At this point it turns out, that in the non-commutative case, a literal free analog of the classical Milstein Method for commutative SDE's with exactly evaluated iterated integrals seems to be impossible. Therefore  in case of nonlinear diffusion $b$ and a general von Neumann algebra $\mathcal{A}$ it remains to apply the so called subdivision method, which we develop in the following. In \Cref{subsec:dgt1nonlinsubdiv} this method will be applied to the general case $d>1$.

\begin{defn}\label{def:subinterval} Let $T>0$. Consider a discretization of $[0,T]$ as in \Cref{def:freeEM-Definition} with $\Delta t<1$.  Set $n=\lceil\frac{1}{\Delta t}\rceil$ and $\delta t=\frac{1}{n}\Delta t$. By a subdivision of an interval $[t,t+\Delta t]$ for some $t>0$ we understand the partition $t=\tau_0<\tau_1<\dots<\tau_{n}=t+\Delta t$ of $[t,t+\Delta t]$ into $n\in\N$ intervals.
\end{defn}
In the following we use the abbreviation $\Delta W_{t,\tau}=W_\tau - W_t, \, 0\leq t\leq \tau$. If $t$ and $\tau$ are discretization points, then we only write the index of the discretization points as indices of $\Delta W$ (same for $X_t$).  Consider a subdivision of $[t,t+\Delta t]$ as described in \Cref{def:subinterval}. We use $\Delta W_{l-1,l}=W_{\tau_{l}}-W_{\tau_{l-1}}$ frequently. If one of the discretization points coincide with the end points of the interval ($\tau_0=t$), we write the point in the index, e.g. $\Delta W_{t,l-1}=W_{\tau_{l-1}}-W_t$.
\begin{thm}\label{def:sudivdequal1} Let $d=1$. 
	As in \Cref{def:freeEM-Definition} consider a partition of $[0,T]$ into $L\in\N$ intervals $[t_{k},t_{k+1}],k=0,\dots,L-1$ with constant step size $\Delta t=\frac{T}{L}$.  Consider a subdivision of the interval $[t_k,t_{k+1}]$ as defined in \Cref{def:subinterval}. The approximation is as follows.
	\begin{equation}
		\oX_{k+1}=\oX_{k}+a(\oX_{k})\Delta t+ b(\oX_{k})\Delta W_{{k}}b(\oX_{k}) + \om_k(\Delta t),
	\end{equation}
	where $\om(\Delta t)$ simulates the iterated Ito-integrals, i.e. 	
	\begin{align*}
		\om_k(\Delta t) &= b(\oX_k)\left\{
		\sum\limits_{l=1}^n  \left[b(\oX_k+\overline{V}_l)-b(\oX_k) \right]\right\}
		b(\oX_k)\\
		&+\sum\limits_{l=1}^n  b(\oX_k)\left\{\Delta W_{l-1,l}\left[b(\oX_k+ b(\oX_k)\Delta W_{t,l-1})-b(\oX_k)\right]\right\}b(\oX_k)\\
		&+\sum\limits_{l=1}^n b(\oX_k)\left\{\left[b(\oX_k+\Delta W_{t,l-1}b(\oX_k))-b(\oX_k)\right]\Delta W_{l-1,l}\right\}b(\oX_k),
	\end{align*}
	and $\overline{V}_l=\Delta W_{l-1,l}b(\oX_k)-\varphi(b(\oX_k))\Delta t$. Then this method shows strong convergence with order $\gamma=1$.
\end{thm}
\begin{proof} Due to \Cref{thm:thm-order1-general} it suffices to show, that $\om_k(\Delta t)$ has the property \eqref{E:Estim-of-m-order-1}. This follows from \Cref{lem:diffmtermsproofthm} by setting $d=1$.
\end{proof}

\subsection{\texorpdfstring{Nonlinear Diffusion, $d>1$}{Nonlinear Diffusion, d}}\label{subsec:dgt1nonlinsubdiv}

The key in the development of the methods in \Cref{def:fMMlinearDiff} and \Cref{thm:deq1-matrix} was the application of the It\^{o} formula to resolve the iterated integrals to a product to obtain \Cref{E:fMMvariante1-def-2x}. In the general case, we have to deal with
\begin{multline}
	m_t^i(\Delta t) = \int_{t}^{t+\Delta t} b_t^i dW_s\left( 	T_{c^{i,[1]}}^{X_t,X_t}\left(\sum\limits_{j=1}^{d}b_t^j \int_{t}^s dW_uc_t^j \right)\right) + \\ + 
	\int_{t}^{t+\Delta t}\left(T_{b^{i,[1]}}^{X_t,X_t}\left(\sum\limits_{j=1}^{d}b_t^j\int_{t}^sdW_uc_t^j\right)\right)dW_s c_t^i.
\end{multline}
To simplify the iterated integrals by the It\^{o} formula in product form, it is algebraically required that $c^i=b^j, i,j=1,\dots d$, but this is equivalent to $d=1$.  Unfortunately, the application of It\^{o} is not possible for $d>1$.\\
The general case $d>1$ requires therefore the treatment of the iterated integrals. So far, to the best knowledge of the authors, there is no approximation method of non-commutative iterated free stochastic integrals known. Therefore we treat the integrals by the subdivision method developed in the previous chapter.

\subsubsection{\texorpdfstring{Subdivision Method, $d>1$}{Subdivision Method, d}}\label{sec:defoffSMandThm}
In this section we define a method of strong order $\gamma=1$ for general $\mathcal{A}$ and $d>1$. It is based on approximating the iterated integrals in \eqref{E:m1rep1-1-def} by a proper subdivision of $\Delta t$, such that, loosely speaking, the iterated integrals are approximated good enough to obtain a higher convergence rate of the numerical method.\\
As a first step, we simplify \eqref{E:m1rep1-1-def}. 
Since elements $b^i(X_t), c^i(X_t)$ and $e^{is X_t}$ belong to the same subalgebra generated by $X_t\in\mathcal{A}^{sa}$, they do commute and due to the linearity of the operator integrals. We can rearrange \eqref{E:m1rep1-1-def} to
\begin{multline}\label{E:m1rep1-1-def-1}
	\sum\limits_{i=1}^d m_t^i(\Delta t) = \sum\limits_{i=1}^d \sum\limits_{j=1}^d b_t^i\int_t^{t+\Delta t}dW_s T_{c^{i,[1]}}^{X_t,X_t}(b_t^j\int_t^sdW_u)c_t^j
	+\\+
	\sum\limits_{i=1}^d\sum\limits_{j=1}^d b_t^j \int_t^{t+\Delta t}T_{b^{i,[1]}}^{X_t,X_t}(\int_t^sdW_u c_t^j)dW_sc_t^i
	=\\=
	\sum\limits_{i=1}^d \sum\limits_{j=1}^d b_t^i\left( \int_t^{t+\Delta t}dW_sT_{c^{i,[1]}}^{X_t,X_t}(b_t^j\int_t^s dW_u) \right. +\\+ \left. \int_t^{t+\Delta t}T_{b^{j,[1]}}^{X_t,X_t}(\int_t^sdW_uc_t^i) dW_s\right)c_t^j
	=\\=
	\sum\limits_{i=1}^d \sum\limits_{j=1}^d b_t^i\left(I_1^{i,j}(\Delta t) + I_2^{i,j}(\Delta t)\right)c_t^j
\end{multline}
Now consider a subdivision of $[t,t+\Delta t]$ as described in \Cref{def:subinterval}. We rewrite $I_1^{i,j}$ and $I_2^{i,j}$ using the step-width $\delta t \leq \Delta t^2$, 
\begin{align*}
	I_1^{i,j}(\Delta t)
	&=
	\sum_{l=1}^n\int_{\tau_{l-1}}^{\tau_l}dW_sT_{c^{i,[1]}}^{X_t,X_t}(b_t^j\int_t^s dW_u)\\
	&=
	\sum_{l=1}^n\int_{\tau_{l-1}}^{\tau_l}dW_sT_{c^{i,[1]}}^{X_t,X_t}(b_t^j\int_{\tau_{l-1}}^s dW_u + b_t^j\int_t^{\tau_{l-1}}dW_u)\\
	&=\sum_{l=1}^n\left[\int_{\tau_{l-1}}^{\tau_l}dW_sT_{c^{i,[1]}}^{X_t,X_t}(b_t^j\int_{\tau_{l-1}}^s dW_u)+\Delta W_{l-1,l}T_{c^{i,[1]}}^{X_t,X_t}(b_t^j\Delta W_{t,l-1})\right] 
\end{align*}

Due to freeness we have the following estimation by applying \cite[Theorem 4.3.8]{Skripka2019} and considering that $\left\|\Delta W_{t,t+\delta t} \right\|^2=\mathcal{O}(\delta t^2)$,
\begin{multline}\label{eq:L2estimreferee}
	\|\int_{t}^{t+\delta t}T^{X_t,X_t,X_t}_{f^{[2]}}(dW_sX_t-X_tdW_s,b_t\int_t^sdW_u)\|^2
	=\\=\|\int_{t}^{t+\delta t}T^{X_t,X_t,X_t}_{f^{[2]}}(\int_t^sdW_ub_t,X_tdW_s-dW_sX_t)\|^2 =O(\delta t^2).
\end{multline}
Now we are ready to apply \Cref{lem:referee} to the first summand of $I_1^{i,j}(\Delta t)$. By \Cref{lem:referee} we push the integration variable $dW_s$ into the operator integral $T_{c^{i,[1]}}^{X_t,X_t}(\cdot)$ in the first summand.  Together with \eqref{eq:L2estimreferee} we obtain
\begin{multline*}
	I_1^{i,j}(\Delta t)
	=\sum_{l=1}^n\left[T_{c^{i,[1]}}^{X_t,X_t}(\int_{\tau_{l-1}}^{\tau_l}dW_sb_t^j\int_{\tau_{l-1}}^s dW_u)
	\right.+\\+\left.\Delta W_{l-1,l}T_{c^{i,[1]}}^{X_t,X_t}(b_t^j\Delta W_{t,l-1})\right]
	+A_1^{i,j}(\delta t).
\end{multline*}
The term $A_1^{i,j}(\delta t)$ is the sum over $l=1,\dots,n$ of the triple operator integral produced through the application of \Cref{lem:referee}.  
Due to \eqref{eq:L2estimreferee} it follows $\|A_1^{i,j}(\delta t)\|^2 =  \mathcal{O}(n\delta t^2)=\mathcal{O}(\Delta t^3)$, where the constant is independent of the index $l$. Similarly we obtain 
\begin{multline*}
	I_2^{i,j}(\Delta t)
	=
	\sum_{l=1}^n\left[T_{b^{j,[1]}}^{X_t,X_t}(\int_{\tau_{l-1}}^{\tau_l}\int_{\tau_{l-1}}^s dW_uc_t^idW_s)
	\right. +\\+\left. T_{b^{j,[1]}}^{X_t,X_t}(\Delta W_{t,l-1}c_t^i)\Delta W_{l-1,l}\right]+A_2^{i,j}(\delta t).
\end{multline*}
As for $A_1^{i,j}(\delta t)$ we conclude $\|A_2^{i,j}(\delta t)\|^2=\mathcal{O}(\Delta t^3)$.
Finally \eqref{E:m1rep1-1-def-1} is rewritten as
\begin{equation}\label{eq:mTermStochTaylor}
	\begin{aligned}
	m_t^i(\Delta t) &= \sum\limits_{j=1}^d b_t^i\left(I_1^{i,j}(\Delta t) + I_2^{i,j}(\Delta t)\right)c_t^j =\\
	&= 
	\sum\limits_{j=1}^d \sum\limits_{l=1}^n b_t^i\left[
	T_{c^{i,[1]}}^{X_t,X_t}(\int_{\tau_{l-1}}^{\tau_{l}}dW_sb_t^j\int_{\tau_{l-1}}^s dW_u)\right]c_t^j + \\
	&+
	\sum\limits_{j=1}^d \sum\limits_{l=1}^n b_t^i\left[T_{b^{j,[1]}}^{X_t,X_t}(\int_{\tau_{l-1}}^{\tau_l}\int_{\tau_{l-1}}^s dW_uc_t^idW_s)
	\right]c_t^j + \\
	&+\sum\limits_{j=1}^d \sum\limits_{l=1}^n b_t^i\left[  \Delta W_{l-1,l}T_{c^{i,[1]}}^{X_t,X_t}(b_t^j\Delta W_{t,l-1}) \right]c_t^j +\\
	&+ \sum\limits_{j=1}^d \sum\limits_{l=1}^n b_t^i\left[T_{b^{j,[1]}}^{X_t,X_t}(\Delta W_{t,l-1}c_t^i)\Delta W_{l-1,l} \right]c_t^j\\
	&+ \sum\limits_{j=1}^d (A_1^{i,j}(\delta t)+A_2^{i,j}(\delta t)).
\end{aligned}
\end{equation}
In the general case $d>1$ we have to deal with the condition $b_t^j\neq c_t^i$. The consequence is, that the sum of the operator integrals $T_{c^{i,[1]}}^{X_t,X_t}(\cdot)$ and $T_{b^{j,[1]}}^{X_t,X_t}(\cdot)$ in \eqref{eq:mTermStochTaylor} cannot be combined. Only the case $d=1$ allows this sum to be simplified. Then, a simplification of iterated stochastic integrals by the It\^{o} formula is possible. For details, we refer to \Cref{sec:dequal1} but for now, our aim is to define a numerical method suitable for the case $d>1$. To do so, we now discretize \eqref{eq:mTermStochTaylor} by \Cref{lem:discrete-opt-int} and skip the terms $A_1^{i,j}$ and $A_2^{i,j}$ in \eqref{eq:mTermStochTaylor}.
\begin{defn}[fSM]\label{def:sudivgeneral}
	As in \Cref{def:freeEM-Definition} consider a partition of $[0,T]$ into $L\in\N$ intervals $[t_{k},t_{k+1}],k=0,\dots,L-1$ with constant step size $\Delta t=\frac{T}{L}$. Consider a subdivision of $[t_k,t_{k+1}]$ as defined in \Cref{def:subinterval}. If we use the abbreviation $b^i(X_k)=b_k^i,\, b^i(\oX_k)=\ob_k^i$ (same for $a, c^i$), then the  fSM approximation of the solution of $X_t$ is defined as
	\begin{equation}\label{eq:fSMMethode}
		\oX_{k+1}=\oX_{k}+\oa_k\Delta t+ \ob_k\Delta W_{{k}}\ob_k + \sum\limits_{i=1}^d \om_k^i(\Delta t),
	\end{equation}
	where $\om^i(\Delta t)$ simulates the iterated Ito-integrals,
	\begin{equation}\label{eq:fSM-m-term}
	\begin{aligned}
		\om_k^i(\Delta t) &= \sum\limits_{j=1}^d \sum\limits_{l=1}^n \ob^i_k
		  \left[c^i\left(\oX_k+\Delta W_{l-1,l} \ob^j_k\Delta W_{l-1,l}\right)-\oc^i_k\right]\oc^j_k + \\
		&+ \sum\limits_{j=1}^d \sum\limits_{l=1}^n \ob^i_k
		\left[b^j\left(\oX_k + \Delta W_{l-1,l}\oc^i_k\Delta W_{l-1,l}\right)-\ob^j_k\right]\oc^j_k\\
		&+\sum\limits_{i=1}^d \sum\limits_{l=1}^n  \ob^i_k\left\{\Delta W_{l-1,l}\left[c^i(\oX_k+ \ob^j_k\Delta W_{t_k,l-1})-\oc^i_k\right]\right\}\oc^j_k\\
		&+\sum\limits_{i=1}^d \sum\limits_{l=1}^n \ob^i_k\left\{(b^j(\oX_k+\Delta W_{t_k,l-1}\oc^i_k)-\ob^j_k)\Delta W_{l-1,l}\right\}\oc^j_k
	\end{aligned}
	\end{equation}
\end{defn}

\begin{lem}\label{lem:diffmtermsproofthm} Let $a,b^i, c^i\in W_3(\R)$.  Consider the terms
	\eqref{eq:mTermStochTaylor} and \eqref{eq:fSM-m-term}.
	Assuming a bounded numerical approximation $\oX_k$, i.e. there is a  
	$M>0$, such that $\|\oX_k\|<M<\infty$ for all $L\in\N, L>1$ and $k=0,\dots, L-1$, then there is an estimation of the form
	\begin{equation}\label{eq:propfSMdiffmterms}
		\|m_k^i(\Delta t)-\om_k^i(\Delta t)\|^2\leq K_1 \Delta t \|X_k-\oX_k\|^2 + K_2\Delta t^3,
	\end{equation}
	where the constants $K_1,K_2>0$ are independent of the discretization $\Delta t$ resp. $\delta t$.
\end{lem}
\begin{proof} Consider the time point $t=t_k$. We use the abbreviations $X_k=X_{t_k}$, $b^i(\oX_k)=\ob^i_k$, $b^i(X_k)=b^i_k$. We will make use of the inequality $(\sum_{k=1}^m a_k)^2\leq m \sum_{k=1}^m a_k^2, m\in\N,\, a_k\in\R$ several times.
	Then \eqref{eq:mTermStochTaylor},  evaluated at $t=t_k$ reads
	\begin{equation}\label{mthermhelp}
		\begin{aligned}
			m_k^i(\Delta t) &=
			\sum\limits_{j=1}^d \sum\limits_{l=1}^n b_k^i\left[
			T_{c^{i,[1]}}^{X_k,X_k}(\int_{\tau_{l-1}}^{\tau_l}dW_sb_k^j\int_{\tau_{l-1}}^s dW_u)\right]c_k^j+\\ 
			&+\sum\limits_{j=1}^d \sum\limits_{l=1}^n b_k^i\left[
			T_{b^{j,[1]}}^{X_k,X_k}(\int_{\tau_{l-1}}^{\tau_l}\int_{\tau_{l-1}}^s dW_uc_k^idW_s)
			\right]c_k^j\\
			&+\sum\limits_{j=1}^d \sum\limits_{l=1}^n b_k^i\left[  \Delta W_{l-1,l}T_{c^{i,[1]}}^{X_k,X_k}(b_k^j\Delta W_{t_k,l-1})\right]c_k^j + \\
			&+\sum\limits_{j=1}^d \sum\limits_{l=1}^n b_k^i\left[
			T_{b^{j,[1]}}^{X_k,X_k}(\Delta W_{t_k,l-1}c_t^i)\Delta W_{l-1,l} \right]c_k^j\\
			&+ \sum\limits_{j=1}^d (A_1^{i,j}(\delta t)+A_2^{i,j}(\delta t)).
		\end{aligned}
	\end{equation}
	For completeness, we repeat \eqref{eq:fSM-m-term}:
	\begin{equation}\label{eq:omtermhelp}
		\begin{aligned}
			\om_k^i(\Delta t) &= \sum\limits_{j=1}^d \sum\limits_{l=1}^n \ob^i_k
			\left[c^i\left(\oX_k+\Delta W_{l-1,l} \ob^j_k\Delta W_{l-1,l}\right)-\oc^i_k\right]\oc^j_k + \\
			&+ \sum\limits_{j=1}^d \sum\limits_{l=1}^n \ob^i_k
			\left[b^j\left(\oX_k + \Delta W_{l-1,l}\oc^i_k\Delta W_{l-1,l}\right)-\ob^j_k\right]\oc^j_k\\
			&+\sum\limits_{i=1}^d \sum\limits_{l=1}^n  \ob^i_k\left\{\Delta W_{l-1,l}\left[c^i(\oX_k+ \ob^j_k\Delta W_{t_k,l-1})-\oc^i_k\right]\right\}\oc^j_k\\
			&+\sum\limits_{i=1}^d \sum\limits_{l=1}^n \ob^i_k\left\{(b^j(\oX_k+\Delta W_{t_k,l-1}\oc^i_k)-\ob^j_k)\Delta W_{l-1,l}\right\}\oc^j_k
		\end{aligned}
	\end{equation}
	By the help of \Cref{lem:discrete-opt-int} we can express the difference in the bracket in the last line of \eqref{eq:omtermhelp} in terms of an operator integral. \Cref{lem:discrete-opt-int} states that
	\begin{equation}\label{eq:opdiffhelp}
		T_{c^{i,[1]}}^{\oX_k,\oX_k}(\oU_l) =
		c^{i}(\oX_k+\oU_l)-c^i(\oX_k)+R_{2,c,\oX_k}(\oU_l),
	\end{equation} 
	where $\oU_l=\int_{\tau_{l-1}}^{\tau_l}dW_s\ob_k^j\int_{\tau_{l-1}}^s dW_u$ and $\|R_{2,c,\oX_k}(\oU_l)\|^2=\mathcal{O}(\delta t^4)$. Inserting this  relationship into \eqref{eq:omtermhelp} the right expression in \eqref{eq:opdiffhelp}, after building the sum over all differences, can be estimated due to the freeness of $\oU_l$ to expressions evaluated at $t_k$ as (note the boundedness of the operator functions $b^i,c^i$ in operator norm and boundedness of the numerical solution) 
	\begin{equation}\label{eq:Delta1innorm-rechterTerm}
		\|\sum_i\sum_l \ob  R_{2,c,\oX_k}(\oU_l) \oc\|^2 \leq C_1 d\sum_l\| R_{2,c,\oX_k}(\oU_l)\|^2  
		 \leq
		   C_2 n\delta t^4 \leq C_3 \Delta t^3. 
	\end{equation}
	We will consider each row separately and  go into detail just for the first term in \eqref{mthermhelp} resp. \eqref{eq:omtermhelp}. The other rows in  \eqref{mthermhelp} resp. \eqref{eq:omtermhelp} can be treated similarly. 
	First lets give an estimation of $\Delta_1$, which is the difference of the first summands in the first line of \eqref{mthermhelp} and \eqref{eq:omtermhelp},
	\begin{equation}\label{eq:beweismilsteintermeh1}
		\begin{aligned}
			\Delta_1 = &\sum\limits_{j=1}^d \sum\limits_{l=1}^n b_k^i\left[
			T_{c^{i,[1]}}^{X_k,X_k}(\int_{\tau_{l-1}}^{\tau_l}dW_sb_k^j\int_{\tau_{l-1}}^s dW_u)\right]c_k^j\\
			-&
			\sum\limits_{j=1}^d \sum\limits_{l=1}^n \ob^i_k
			\left[c^i(\oX_k+\Delta W_{l-1,l} \ob^j_k\Delta W_{l-1,l})-\oc^i_k
			\right]
			\oc_{k}^j
			=\\
			=&
			\sum\limits_{j=1}^d \sum\limits_{l=1}^n b_k^i\left[
			T_{c^{i,[1]}}^{X_k,X_k}(\int_{\tau_{l-1}}^{\tau_l}dW_sb_k^j\int_{\tau_{l-1}}^s dW_u)\right]c_k^j\\
			-&
			\sum\limits_{j=1}^d \sum\limits_{l=1}^n \ob^i_k
			\left[
			T_{c^{i,[1]}}^{\oX_k,\oX_k}(\Delta W_{l-1,l} \ob^j_k\Delta W_{l-1,l})-R_{2,c,\oX_k}(\oU_l)\right]\oc_{k}^j 	\\
			=&\sum_{j=1}^d\sum\limits_{l=1}^n \left[b_k^iT_{c^{i,[1]}}^{X_k,X_k}(I_l) c_k^j- \ob_k^iT_{c^{i,[1]}}^{\oX_k,\oX_k}(\oI_l )\oc_k^j \right]	- \sum\limits_{j=1}^d \sum\limits_{l=1}^n\ob^i_k R_{2,c,\oX_k}(\oU_l)\oc^j_k
		\end{aligned}
	\end{equation}
	where $I_l=\int_{\tau_{l-1}}^{\tau_l}dW_sb_k^j\int_{\tau_{l-1}}^s dW_u$,\, $\oI_l=\Delta W_{l-1,l} \ob^j_k\Delta W_{l-1,l}$.\\
	Using \Cref{lem:estim-optint-diff} we can exchange $b_k^iT_{c^{i,[1]}}^{X_k,X_k}(I_l)c_k^j$ by $b_k^iT_{c^{i,[1]}}^{\oX_k,\oX_k}(I_l)c_k^j$ with a penalty term of square operator norm less than $C_1\delta t^2\|X_k-\oX_k\|^2$ 
	(note that the $b_k^i,c_k^i$ are bounded for $k=1,\ldots,L-1$ and $i=1,\ldots, d$). With the technique to express $bdWb-cdWc$ by three symmetric terms (see \Cref{sec:FreeIto}) we obtain the estimation
	\begin{equation*}
		\|b_k^i T_{c^{i,[1]}}^{\oX_k,\oX_k}(I_l) c_k^j-\ob_k^iT_{c^{i,[1]}}^{\oX_k,\oX_k}(I_l)\oc_k^j\|\le C_2 \delta t^2\|X_k-\oX_k\|^2
	\end{equation*}
	(note that $\|T_{c^{i,[1]}}^{\oX_k,\oX_k}(I_l)\|^2\le \text{const}\ \delta t^2$). Thus, we can replace $b_k^i T_{c^{i,[1]}}^{\oX_k,\oX_k}(I_l)c_k^j$ by $\ob_k^iT_{c^{i,[1]}}^{\oX_k,\oX_k}(I_l)\oc_k^j$ with an additional term of square $L_\infty(\varphi)$-norm less or equal to $C_3 \delta t^2 \|X_k-\oX_k\|^2$.
	Hence, with the freeness argument
	\begin{equation}
		\left\|\sum_{l=1}^n\left[b_k^iT_{c^{i,[1]}}^{X_k,X_k}(I_l) c_k^j- \ob_k^iT_{c^{i,[1]}}^{\oX_k,\oX_k}(I_l )\oc_k^j \right]\right\|^2\le C_3\Delta t^3\|X_k-\oX_k\|^2.
	\end{equation}
	Finally,
	$\Delta I_l :=I_l-\oI_l= \int_{\tau_{l-1}}^{\tau_l}dW_s\int_{s}^{l}(b_k^i-\ob_k^i)dW_u+(\oU_l-\oI_l)$, hence $\|\Delta I_l\|^2\le C(\delta t^2\|X_k-\oX_k\|^2+\delta t^2)$ (note that $\|\oU_l-\oI_l\|^2=\mathcal O(\delta t^2)$). The fact that $T_{c^{i,[1]}}^{\oX_k,\oX_k}$ is a bounded linear operator gives us 
	\begin{align*}
		&\|\sum_{l=1}^n\left[b_k^iT_{c^{i,[1]}}^{X_k,X_k}(I_l) c_k^j- \ob_k^iT_{c^{i,[1]}}^{\oX_k,\oX_k}(\oI_l )\oc_k^j \right]\|^2\\
		&\le \|  \sum_{l=1}^n\left[b_k^iT_{c^{i,[1]}}^{X_k,X_k}(I_l) c_k^j- \ob_k^iT_{c^{i,[1]}}^{\oX_k,\oX_k}(I_l )\oc_k^j \right]\|^2 +\|\ob_k^iT_{c^{i,[1]}}^{\oX_k,\oX_k} (\sum_{l=1}^n(I_l-\oI_l ))\oc_k^j\|^2\\
		&+ C_4\Delta t^3\|X_k-\oX_k\|^2\le C_5\Delta t^3\|X_k-\oX_k\|^2+C_6\Delta t^3.
	\end{align*}
	and the proof is done.
\end{proof}

\begin{thm}\label{thm:strong-conv-fsm} Assume $d\geq 1$. 
	Consider a fSDE \eqref{intro-freeSDE-diffform} with the solution $X_t$ on $[0,T]$ (see \Cref{rem:fSDE-sol-theory}). Let $\oX_k$ be a numerical solution calculated by fSM \eqref{eq:fSMMethode} on $[0,T]$ with a discretization $T=L\Delta t$.  For a fixed $k=0,\ldots,L-1$ we choose the subdivision of $[t_k,t_{k+1}]$ as defined in \Cref{def:subinterval}. Let $a,b^i,c^i\in W_3(\R), \, i=1\dots d$. Then the fSM approximation \eqref{eq:fSMMethode} converges strongly with order $\gamma=1$ to the solution $X_k$
	\begin{equation}\label{eq:strong-order-1}
		\|\oX_k-X_k\|_p\leq C\,\Delta t
	\end{equation}
	for all $1\leq p \leq \infty$.
	The constant $C$ is independent of step size $\Delta t$ and subdivision interval $\delta t$.\\
	Furthermore, the numerical solution is uniformly bounded for each value  $k=0,\dots,L-1, \, L\in\N$, i.e. there is a constant $\overline{M}>0$ such that $\|\oX_k\|<\overline{M}$, where $\overline{M}$ does not depend on $L$ resp. $\Delta t$ and the subdivision interval $\delta t$.\\
\end{thm}
\begin{proof}
	Strong convergence follows from \Cref{thm:thm-order1-general} and \Cref{lem:diffmtermsproofthm}. The boundedness of the numerical solution follows with the same arguments as in the proof of \Cref{thm:strong-conv-bound-num-sol}.
\end{proof}

\section{Numerical Examples}\label{sec:numex}
In the following examples we present different fSDEs to confirm the insights from the previous chapters numerically. We select three fSDEs with $d=1$, $d=2$ and different functions $b,c$. Additionally we show that an a posteriori estimation of the convergence order gives the expected results. We mainly follow \cite{SchlueWib2023} and numerically seek an approximation of the solution of the fSDE at a final time $T>0$.
\subsection{A Simple Example}
We consider the fSDE
\begin{equation}\label{ex:explosive-sde}
	dX_t= X_tdW_tX_t
\end{equation}
with start value $X_0=I$. In \cite[Proposition 3.9]{kargin} it is shown, that the spectral distribution of the solution $X_t\in\mathcal{A}^{sa}$ exists of all $t\leq1$ and is supported on the interval 
\begin{equation}
	\left[\frac{(1-\sqrt{t})^2}{(1-t)^2},\frac{(1+\sqrt{t})^2}{(1-t)^2}\right].
\end{equation}
For $t\in ]0,1]$, the density of the spectral distribution is given by
\begin{equation}\label{kap8-example-PDF}
	f(x) = \frac{\sqrt{-(1-T)^2 x^2+2(1+T)x-1}}{2\pi T x^3}.
\end{equation}
For $t=1$ the density is supported on $[1/4,\infty[$.\\
We implement the method \eqref{eq:def-method-linear-d1} with $\lambda=1, \mu=0$ defined in \Cref{def:fMMlinearDiff}. 
The probability density function (PDF) of the eigenvalues of the numerical solution are then an approximation of \eqref{kap8-example-PDF}. 
\begin{figure}[thbp]
	\centering
	\begin{tikzpicture}
		  \pgfplotstableread{data/ExplosiveStrongSimple.dat}\ExplosiveStrongSimple
		  \pgfplotstableread{data/ExplosiveStrongfEMM.dat}\ExplosiveStrongfEMM		
			\begin{loglogaxis}[scale=0.6,
				xlabel = {$\Delta t_R$},
				ylabel = {$\epsilon^{\text{Sim}}$},
				xmin=0.001, xmax=0.05, ymin=0.001, ymax=0.1,
				grid=both,
				minor grid style={gray!25},
				major grid style={gray!25},
				legend entries = {\scriptsize Method $\eqref{eq:def-method-linear-d1}$: $\gamma=0.94$, \scriptsize fEMM: $\gamma=0.49$},
				legend pos = outer north east,
				]
				\addplot table[x=x,y=y1] {\ExplosiveStrongSimple};
				\addplot table[x index=0, y index = 1] {\ExplosiveStrongfEMM};				
			\end{loglogaxis}
		\end{tikzpicture}
	\caption{Strong  convergence properties of fEMM and Method $\eqref{eq:def-method-linear-d1}$ applied to the fSDE \eqref{ex:explosive-sde}. The parameters are $N=10$, $T=0.1$, $S=2570000$.}
	\label{fig:Explosive-Strong}
\end{figure}
The numerical algorithms fEMM and \eqref{eq:def-method-linear-d1} are calculated with four different time steps $\Delta t_R=\frac{T}{2^{13}}2^{R}$ with $R=0,8,9,10$ and final  time $T=0.1$. 
The path of the free Brownian motion in the case $R=0$ will be used to calculate all necessary moments for $R=8,9,10$. The exact solution, which is unkown, will be simulated by $R=0$. An a posteriori error estimation is shown in the succeeding examples. Let's consider $U\in\N$ different paths $W_t(\omega_u), \, u=1,\dots,U$ of the free Brownian motion. In the following $\oX_T(u,\Delta t_R)$ denotes the numerical solution calculated at time point $T$, by time step $\Delta t_R$ and path  $W_t(\omega_u)$. We simply call $W_t(\omega_u)$ as path number $u$.  For each time step $\Delta t_R$ for $R=8,9,10$ and each path number $u$ the strong path-wise error in $L_\infty(\varphi)$-norm is calculated as $e^{\text{Sim}}(u,\Delta t_R)=\|\oX_T(u,\Delta t_R) - \oX_T(u,\Delta t_0)\|$, where the numerical solution $\oX_T(u,\Delta t_R)$ is either calculated by fEMM or the method \eqref{eq:def-method-linear-d1}. Taking the mean value of the path-wise strong error yields
\begin{equation*}
	\epsilon^{\text{Sim}}(\Delta t_R)=1/U\sum_{u=1}^{U}e^{\text{Sim}}(u,\Delta t_R).
\end{equation*}

The comparison of the convergence order of fEMM and method \eqref{eq:def-method-linear-d1} is shown in \Cref{fig:Explosive-Strong}.
\begin{figure}[htpb]
	\centering
	\captionsetup[subfigure]{justification=centering}
	\pgfplotstableread{data/Xthisteins.dat}\Xthisteins
	\pgfplotstableread{data/Xthistanaeins.dat}\Xthistanaeins
	\pgfplotstableread{data/Xthistfuenf.dat}\Xthistfuenf
	\pgfplotstableread{data/Xthistanafuenf.dat}\Xthistanafuenf	
	\subfloat[$t=0.1, N=500$]{
		\begin{tikzpicture}[scale=0.6, xscale=1.05]
			\begin{axis}[ xmin=0.25, xmax=2.5, ymax=2, grid=major] 
				\addplot [ color=black, line width=1pt] table [x index=0, y index=1] {\Xthistanaeins};
				\addplot +[
				hist={density,
					bins=20,		
				}, 
				mark=none,  solid, line width=1pt, color=black, fill=gray, fill opacity=0.3
				] table [y index=0] {\Xthisteins};	
				\addlegendentry{PDF - exact};
			\end{axis}
		\end{tikzpicture}
	}\hspace{0.2cm}
	\subfloat[$t=0.5, N=500$]{
		\begin{tikzpicture}[scale=0.6]
			\begin{axis}[ xmin=0, xmax=6, ymax=2, grid=major]
				\addplot [ color=black, line width=1pt] table {\Xthistanafuenf};
				\addplot +[
				hist={density,
					bins=20,		
				}, 
				mark=none,  solid, line width=1pt, color=black, fill=gray, fill opacity=0.3
				] table [y index=0] {\Xthistfuenf};	
				\addlegendentry{PDF - exact};
			\end{axis}
		\end{tikzpicture}\hspace{0.3cm}
	}
	\caption{Estimation of the probability PDF of the spectral distribution of the exact solution $X_t$ of \eqref{ex:explosive-sde} recovered from it's Cauchy transform at different time points. The bars show an estimation of the PDF via the eigenvalues of the numerical approximation $\oX_t$ calculated by the simple method \eqref{eq:def-method-linear-d1} at different time points. The solid line is the PDF of the exact solution $X_t$.}
	\label{fig:Explosive-hist-t01-t05}
\end{figure}
\Cref{fig:Explosive-hist-t01-t05} shows the density (PDF)  \eqref{kap8-example-PDF} of the spectral distribution of the exact solution $X_t$ of \eqref{ex:explosive-sde} recovered from it's Cauchy transform at two different time points and their approximation calculated by \eqref{eq:def-method-linear-d1} in $\mathcal{M}^{500}(\R)$. The bars show an estimation of the PDF via the eigenvalues of the numerical solution $\oX_t$,  calculated by \eqref{eq:def-method-linear-d1} at different time points. The line is the PDF of the exact solution $X_t$.
\subsection{Geometric Brownian Motion I}
Let's consider the equation
\begin{equation}\label{ex:geo1}
	dX_t= \theta X_tdt + \sqrt{X_t}dW_t\sqrt{X_t}
\end{equation}
with $X_0=I$. This examples demonstrates numerically the order $\gamma=1$ of the method developed in \Cref{thm:deq1-matrix}. In every single time step, the spectral properties of $\oX_k$ has to be computed in order to calculate the single and double operator integrals \eqref{E:Discret-OpInt-Matrix}. The difficulty of this equations is, that if $\theta$ is to small, the eigenvalues turn complex. To confirm $\gamma=1$ we choose $\theta$ large enough. The equations is solved in the von Neumann algebra $\mathcal{M}^5(\R)_{sa}$. \Cref{fig:GeoBrownI-Strong} shows the order of convergence. 

\begin{figure}[htbp]
	\centering
	\begin{tikzpicture}
		\pgfplotstableread{data/GeoIStrongMatrixMethod.dat}\GeoIStrongMatrixMethod
		\pgfplotstableread{data/GeoIStrongSimpleStrongfEMM.dat}\GeoIStrongSimpleStrongfEMM		
		\begin{loglogaxis}[scale=0.6,
			xlabel = {$\Delta t_R$},
			ylabel = {$\epsilon^{\text{Sim}}$},
			xmin=0.0005, xmax=0.005, ymin=0.03, ymax=0.5,
			grid=both,
			minor grid style={gray!25},
			major grid style={gray!25},
			legend entries = {\scriptsize Method $\eqref{eq:deq1-matrix}$: $\gamma=0.99$, \scriptsize fEMM: $\gamma=0.55$},
			legend pos = outer north east,
			]
			\addplot table[x=x,y=y1] {\GeoIStrongMatrixMethod};
			\addplot table[x index=0, y index = 1] {\GeoIStrongSimpleStrongfEMM};					
		\end{loglogaxis}
	\end{tikzpicture}
	\caption{Strong  convergence properties of the method $\eqref{eq:deq1-matrix}$ applied to the fSDE \eqref{ex:geo1}. The parameters are $N=5$, $T=0.1$, $S=500$, $\theta=2$}
	\label{fig:GeoBrownI-Strong}
\end{figure}

The reference solution $\oX_T(u,\Delta t_0)$ was calculated by $\Delta t_0=1/2^{12}$. Further calculations were performed with $\Delta t_R=1/2^{R+7}$ for $R=1,2,3$. We then calculated the error $e^{\text{Sim}}(u,\Delta t_R)=\|\oX_T(u,\Delta t_R) - \oX_T(u,\Delta t_0)\|$ for each path number $u$. As in the previous example, we calculate the mean value of the pathwise error $e^{\text{Sim}}(u,\Delta t_R)$. The plot in \Cref{fig:GeoBrownI-Strong} shows the expected order of convergence.

\subsection{Geometric Brownian Motion II}\label{sec:example-GeoII}
Consider the case $d>1$ with linear diffusion coefficients of the form
\begin{equation}\label{ex:geo2}
	dX_t= \theta X_tdt + X_tdW_t + dW_t X_t
\end{equation}
with start value $X_0=I$.
For analytical insights to the spectral distribution of the solution $X_t$ we again refer to \cite{kargin}.  We apply method fSM \eqref{eq:fSMMethode} and set  $N=10$, $\theta=1$. 
The subdivision method fSM is realized by taking time step values $\Delta t_R = 2^R\cdot 10^{-3}$ for $R=0,1,2,3$ and $\delta t_R=\Delta t_R^2$. We adopt the abbreviations from the first example. The underlying path of the Brownian motion was calculated by a time step $10^{-6}$ ($R=0$). Out of this path, the necessary increments of the Brownian motion for the application of fSM are calculated in dependence of the time step $\Delta t_R$.  As in the previous example, we consider $U\in\N$ different paths. Due to the computational complexity of fSM we perform an a posterior error calculation. For each path $u\in\N$ we calculate
\begin{multline*}
	e^{GeoII}(u,\Delta t_R)=\|\oX_T(u,\Delta t_R)-\oX_T(u,\Delta t_{R-1})\|_1
	=\\=
	\frac{1}{N}\mathbb{E}\left(\text{tr}(|\oX_T(u,\Delta t_R)-\oX_T(u,\Delta t_{R-1})|)\right).
\end{multline*}$\oX_T(u,\Delta t_R)$ denotes the numerical solution by an appropriate method with path $u$ and time step $\Delta t_R$ at time point $T$. We consider $U=2000$ number of different paths and approximate the posterior error as the mean value $ \epsilon^{\text{GeoII}}(\Delta t_R)=\frac{1}{U}\sum_{u=1}^U e^{\text{GeoII}}(u,\Delta t_R)$.
The estimation of the order of convergence is then calculated by
\begin{eqnarray}\label{eq:aposterioriestim}
	\gamma \approx \ogamma(\Delta t_R)=\frac{\log\left(\frac{\epsilon^{\text{GeoII}}(\Delta t_R)}{\epsilon^{\text{GeoII}}(\Delta t_{R-1})}\right)}{\log(2)}.
\end{eqnarray}
The results are listed in \Cref{tab:geo2} and visualized in  \Cref{fig:geo2-strong}. The simulation shows that $\ogamma(\Delta t_3)\approx 0.961$ by fSM and $\ogamma(\Delta t_3)\approx 0.501$ by fEMM.
The results confirm the expected strong convergence orders of $\gamma=0.5$ for fEMM and $\gamma=1$ for fSM.
\begin{figure}[htbp]
	\centering
	\begin{tikzpicture}
		\pgfplotstableread{data/geozweimil.dat}\geozweimil
		\pgfplotstableread{data/geozweimileul.dat}\geozweimileul
		\begin{semilogyaxis}[scale=0.6,
			xlabel = {$R$},
			ylabel = {$\epsilon^{\text{GeoII}}$},
			xmin=0, xmax=4, ymin=0.004, ymax=0.1,
			grid=both,
			minor grid style={gray!25},
			major grid style={gray!25},
			legend entries = {\scriptsize fSM: $\ogamma(\Delta t_3)\approx 0.961$, \scriptsize fEMM: $\ogamma(\Delta t_3)\approx 0.501$},
			legend pos = outer north east,
			]
			\addplot table[x index=0, y index = 1] {\geozweimileul};	
			\addplot table[x index=0, y index = 1] {\geozweimil};
		\end{semilogyaxis}
	\end{tikzpicture}
	\caption{Strong convergence properties of fEMM and fSM applied to the fSDE \eqref{ex:geo2}. Parameters are $N=10$, $T=0.2$, $U=20000$.}
	\label{fig:geo2-strong}
\end{figure}
~\\~
\begin{center}
	\begin{table}[ht]
		\centering
		\begin{tabular}{||c|c|c||} 
			\hline
			$R$ & $\epsilon^{\text{GeoII}}$ by fEMM & $\epsilon^{\text{GeoII}}$ by fSM \\ [0.5ex] 
			\hline\hline
			1 & 3.093e-02 & 8.720e-03 \\ 
			\hline
			2 & 4.376e-02 & 1.688e-02  \\
			\hline
			3 & 6.189e-02 & 3.551e-02 \\
			\hline
		\end{tabular}
		\caption{Data visualized in \Cref{fig:geo2-strong}.}
		\label{tab:geo2}
	\end{table}
\end{center}

\subsection{CIR Equation}
Consider the following nonlinear case with $d>1$,
\begin{equation}\label{eq:CIR}
	dX_t= (a-bX_t)dt + \frac{\sigma}{2}\sqrt{X_t}dW_t + dW_t\frac{\sigma}{2}\sqrt{X_t}
\end{equation}
with start value $X_0=I$. 
For analytical insights to the spectral distribution and the existence of the solution $X_t$ we refer to \cite{freeCIR}. We perform the simulation with $N=2$, $a=I$, $b=0.1$, $\sigma=0.2$. 
Similar to the example in \Cref{sec:example-GeoII}, we proceed by applying method fSM \eqref{eq:fSMMethode} with time steps $\Delta t_R = 2^R\cdot 10^{-3}$ for $R=0,1,2,3$ and $\delta t_R=\Delta t_R^2$. The underlying path of the Brownian motion was calculated by a time step $10^{-6}$. The error for path $u\in\N$ is calculated at
\begin{multline*}
	e^{\text{CIR}}(u,\Delta t_R)=\|\oX_T(u,\Delta t_R)-\oX_T(u,\Delta t_{R-1})\|_1
	=\\=
	\frac{1}{N}\mathbb{E}\left(\text{tr}(|\oX_T(u,\Delta t_R)-\oX_T(u,\Delta t_{R-1})|)\right).
\end{multline*}
 $\oX_T(u,\Delta t_R)$ denotes the numerical solution by fEMM or fSM with path $u$ and time step $\Delta t_R$ at time point $T$. We consider $U=2000$ number of different paths and approximate the a posteriori error as $ \epsilon^{\text{CIR}}(\Delta t_R)=\frac{1}{U}\sum_{u=1}^Ue^{\text{CIR}}(u,\Delta t_R)$.
The convergence order is numerically estimated by
\begin{eqnarray}\label{eq:aposterioriestimCIR}
	\gamma \approx \ogamma(\Delta t_R)=\frac{\log\left(\frac{\epsilon ^{\text{CIR}}(\Delta t_R)}{\epsilon^{\text{CIR}}(\Delta t_{R-1})}\right)}{\log(2)}.
\end{eqnarray}
The results are listed in \Cref{tab:CIR} and visualized in  \Cref{fig:CIR-strong}. The simulation show that $\ogamma(\Delta t_3)\approx 0.999$ by fSM and $\ogamma(\Delta t_3)\approx 0.523$ by fEMM.
\begin{figure}[htbp]
	\centering
	\begin{tikzpicture}
		\pgfplotstableread{data/cirzweimil.dat}\cirzweimil
		\pgfplotstableread{data/cirzweimileul.dat}\cirzweimileul
		\begin{semilogyaxis}[scale=0.6,
			xlabel = {$R$},
			ylabel = {$\epsilon^{\text{CIR}}$},
			xmin=0, xmax=4, ymin=0.00001, ymax=0.001,
			grid=both,
			minor grid style={gray!25},
			major grid style={gray!25},
			legend entries = {\scriptsize fSM: $\ogamma(\Delta t_3)\approx 0.999$, \scriptsize fEMM: $\ogamma(\Delta t_3)\approx 0.523$},
			legend pos = outer north east,
			]
			\addplot table[x index=0, y index = 1] {\cirzweimileul};	
			\addplot table[x index=0, y index = 1] {\cirzweimil};
		\end{semilogyaxis}
	\end{tikzpicture}
	\caption{Strong convergence properties of fEMM and fSM applied to the fSDE \eqref{eq:CIR}. Parameters are $N=2$, $T=0.2$, $U=2000$.}
	\label{fig:CIR-strong}
\end{figure}
~\\~
\begin{center}
	\begin{table}[ht]
		\centering
		\begin{tabular}{||c|c|c||} 
			\hline
			$R$ & $\epsilon^{\text{CIR}}$ by fEMM & $\epsilon^{\text{CIR}}$ by fSM \\ [0.5ex] 
			\hline\hline
			1 & 2.136e-04 & 3.245e-5 \\ 
			\hline
			2 & 3.069e-04 & 6.486e-5  \\
			\hline
			3 & 4.298e-04 & 1.290e-4 \\
			\hline
		\end{tabular}
		\caption{Data visualized in \Cref{fig:CIR-strong}.}
		\label{tab:CIR}
	\end{table}
\end{center}

\bibliographystyle{plainurl}
\bibliography{biblio}

\end{document}